\documentclass[12pt]{article}
\usepackage{epic,amssymb}

\advance \textheight by -1 \topmargin
\advance \textheight by -1 \headheight
\advance \textheight by -1 \headsep
\advance \textheight by -2 \footskip
\vsize = \textheight
\hsize = \textwidth
\newcommand{\vol}{\mbox{\rm vol}}
\newcommand{\rank}{\mbox{\rm rank}}
\newcommand{\mass}{\mbox{\rm mass}}
\newcommand{\Min}{\mbox{\rm Min}}

\newcommand{\Aut}{\mbox{\rm Aut}}

\newtheorem{theorem}{Theorem}[section]

\newtheorem{satz}[theorem]{Theorem}
\newtheorem{prop}[theorem]{Proposition}

\newtheorem{kor}[theorem]{Corollary}
\newtheorem{defi}[theorem]{Definition}

\newcommand{\bew}{\noindent\underline{Proof.}\ }
\newtheorem{bem}[theorem]{Remark}
\newtheorem{rem}[theorem]{Remark}
\newtheorem{lemma}[theorem]{Lemma}

\renewcommand{\setminus}{-}
\newcommand{\fa}{\mbox{ for}\mbox{ all }}
\newcommand{\be}{\begin{enumerate}}
\newcommand{\ee}{\end{enumerate}}
\newcommand{\bi}{\begin{itemize}}
\newcommand{\ei}{\end{itemize}}

\newcommand{\ba}{\begin{array}}
\newcommand{\ea}{\end{array}}

\newcommand{\disj}{\stackrel{.}{\cup}}
\newcommand{\Z}{{\mathbb{Z}}}
\newcommand{\Q}{{\mathbb{Q}}}
\newcommand{\F}{{\mathbb{F}}}
\newcommand{\N}{{\mathbb{N}}}
\newcommand{\R}{{\mathbb{R}}}

\newcommand{\eb}{\phantom{zzz}\hfill{$\square $}\smallskip}

\renewcommand{\em}{\sf}

\newcommand{\knubbel}{\begin{picture}(0,0)(5,-3)  \put(1,1){\circle*{5}} \end{picture}}

\begin{document}

\Large
\begin{center}
{\bf Low dimensional strongly perfect lattices. \\
III: Dual strongly perfect lattices of dimension 14.}
\end{center}
\normalsize
\begin{center}
Gabriele Nebe at Aachen
 and Boris Venkov at St. Petersburg 
\end{center}

\small
{\sc Abstract}: 
The extremal 3-modular lattice $[\pm G_2(3)]_{14}$ with automorphism 
group $C_2 \times G_2(\F_3) $  is the unique dual strongly perfect lattice 
of dimension 14.

\normalsize

\section{Introduction.}

This paper continues the classification of strongly perfect lattices
in \cite{Venkov}, \cite{dim10}, \cite{dim12}, \cite{dim13}. 
A lattice $L$ in Euclidean space $(\R ^n,(,))$ is called 
{\em strongly perfect}, if the set of its minimal vectors forms a 
spherical 4-design. 
The most important property of such lattices is that they provide nice 
examples of {\em extreme} lattices, where
the  density of the associated lattice sphere packing attains 
a local maximum.
In fact \cite{Schuermann} shows that the density of a 
strongly perfect lattice  is even a local maximum on the space of 
all periodic packings, so these lattice packings are periodic extreme. 

All strongly perfect lattices are known up to dimension 12. 
They all share the property of being 
{\em dual strongly perfect}, which means that both lattices 
$L$ and its dual lattice $L^*$ are strongly perfect 
(see Definition \ref{stperf}). 
The only known strongly perfect lattice for which the dual is 
not strongly perfect is $K_{21}'$ (see \cite[Tableau 19.2]{Venkov}) 
in dimension 21. 

One method to show that a lattice $L$ is strongly perfect is to
use its automorphism group $G=\Aut(L)$, the stabilizer of $L$ in
the orthogonal group. If this group has no harmonic 
invariant of degree $\leq 4$, then all  $G$-orbits are 
spherical 4-designs, and so is 
$$\Min (L) := \{ \ell \in L \mid (\ell,\ell ) = \min (L) \} .$$
Since $\Aut(L) = \Aut(L^*) $ such lattices are also 
dual strongly perfect. 
Moreover all non empty layers 
$$L_a := \{ \ell \in L \mid (\ell , \ell ) = a \} $$ 
are 4-designs, since they are also disjoint unions of $G$-orbits.
We call such lattices {\em universally perfect} (see Definition \ref{stperf}).
Though being
 universally perfect involves infinitely many layers of the lattice 
it can be checked with a finite computation using the theory of modular 
forms. 
For any harmonic polynomial $p$ of degree $d$ and a lattice $L\subset(\R^n,(,))$
 the theta series 
$$\theta _{L,p} := \sum _{\ell \in L} p(\ell )  q^{(\ell,\ell )} $$
is a modular form of weight $k:=\frac{n}{2} +d $ and hence lies in a 
finite dimensional vector space $M_k(\ell )$ (depending on the level $\ell$
 of $L$).
A lattice is universally perfect, if and only if this theta series 
$\theta _{L,p}$ 
vanishes for all harmonic polynomials $p$ of degree 2 and 4, a condition
that can be tested from the first few Fourier coefficients of 
$\theta _{L,p}$.

Under certain conditions, the theory of modular forms allows to show 
that $\theta _{L,p} = 0$ for all harmonic polynomials 
$p$ of degree 2 and 4 which again shows that all layers of $L$ are 
spherical 4-designs (see for instance \cite{BachocVenkov}, \cite[Theor\`eme 16.4]{Venkov}). 
In fact this is the only known way to conclude that all 
even unimodular lattices of dimension 32 and minimum 4 provide 
locally densest lattice sphere packings
see \cite[Section 16]{Venkov}. 
By \cite{king} there are more than $10^6$ such lattices and a complete
classification is unknown. 

If $L$ is a universally perfect lattice, then 
 the theta-transformation formula shows that also 
$\theta _{L^*,p} = 0 $ for all harmonic polynomials of degree 2 and 4
and hence $L$ is dual strongly perfect. 
So a classification of all dual strongly perfect lattices includes 
those that are universally perfect. 

Universally perfect lattices also play a role in Riemannian geometry. 
If $L$ is a universally perfect lattice (of fixed covolume 
$\vol (\R^n/L) ^2 = \det(L) = 1$, say) then 
the torus
$\R^n/L^*$ defined by the dual lattice $L^*$
provides a strict local minimum of the
height function on the set of all $n$-dimensional flat tori of
volume 1 (\cite[Theorem 1.2]{coulange1}).
R. Coulangeon also shows that universally perfect lattices $L$
achieve local minima of Epstein's zeta function,
they are so called $\zeta $-extreme lattices. 
The question to find $\zeta $-extreme lattices has a 
long history going back to Sobolev's work \cite{Sobolev}
on numerical integration and \cite{Delone}.

Restricting to dual strongly perfect lattices gives us additional
powerful means for classification and non-existence proofs:
Our general method to classify all strongly perfect lattices in 
a given dimension usually starts with a finite list of possible 
pairs $(s,\gamma )$, where $s = s(L) = \frac{1}{2} |\Min (L) |$ is half 
of the kissing number of $L$ and 
$$\gamma  = \gamma '(L)^2 = \gamma '(L^*)^2 = \min (L) \min (L^*)$$
 the Berg\'e-Martinet invariant of $L$.
For both quantities there are good upper bounds given 
in \cite{Anstreicher} resp. \cite{Elkies}. 
Note that $\gamma $ is just the product of the values of the
Hermite function on $L$ and $L^*$.
Using the general equations for designs given in Section \ref{DDD} 
a case by case analysis allows either to exclude certain of the
possibilities $(s,\gamma )$ or to factor $\gamma = m \cdot r $ such that 
rescaled to minimum $\min (L^*) = m$, the lattice $L^*$ is integral 
(or even) and in particular contained in its dual lattice $L$ (which is 
then of minimum $r$). 
For dual strongly perfect lattices we can use a similar argumentation 
to obtain a finite list of possibilities $(s',\gamma )$ for 
$s' = s(L^*)$ and in each case a factorization $\gamma = m' \cdot r' $ 
such that $L$ is integral (or even) if rescaled to $\min (L) = m' $. 
But this allows to obtain the exponent (in the latter scaling) 
$$\exp (L^*/L)  = \frac{m}{r'} $$
which either allows a direct classification of all such lattices 
$L$ or the use  of modular forms to exclude 
the existence of a modular form 
$\theta _{L}=\theta _{L,1} $ of level $\frac{m}{r'} $ and weight $\frac{n}{2}$ 
starting with 
$1+2sq^{m'} + \ldots $, such that its image under the Fricke involution starts 
with $1+2s'q^m + \ldots $ and both $q$-expansions have nonnegative 
integral coefficients. 

{\bf Acknowledgement.} 
Most of the research was done during Boris Venkov's visit to 
the RWTH Aachen university, from October 2007 to July 2008. 
We thank the Humboldt foundation for their financial support.

\section{Some general equations.}{\label{allg}}

\subsection{General notation.}{\label{not}}

For a lattice $\Lambda $ in $n$-dimensional Euclidean space $(\R^n , (,))$
we denote by
$$\Lambda ^* := \{ v \in \R^n \mid (v,\lambda) \in \Z \mbox{ for all } 
\lambda\in\Lambda \} $$ its dual lattice and 
by
$$\Lambda _a :=\{ \lambda \in \Lambda \mid (\lambda , \lambda ) = a \} $$
the vectors of square length $a$. 
In particular, if 
$\min := \min (\Lambda ) := \min \{ (\lambda ,\lambda ) \mid 
0\neq \lambda \in \Lambda \}$ denotes the {\em minimum} of the lattice then
$\Lambda _{\min } = \Min (\Lambda )$
 is the set of minimal vectors in $\Lambda $, its
cardinality is known as the {\em kissing number} $2s(\Lambda )$ 
of the lattice $\Lambda$.
There are general bounds on the 
 kissing number, of an $n$-dimensional lattice. 
For $n=14$ the bound is
 $ |\Lambda _{\min } |\leq 2\cdot 1746$
and hence $s(\Lambda ) \leq 1746 $
 (see \cite{Anstreicher}).

Let $$\gamma _n := \max \{ \frac{\min(\Lambda )}{\det(\Lambda )^{1/n}} \mid
\Lambda \mbox{ is an $n$-dimensional lattice } \} $$
denote the {\em Hermite constant}.
The precise value for $\gamma _n$ is known for $n\leq 8$ and $n=24$.
However upper bounds are given in \cite{Elkies}. 
In particular $\gamma _{14} \leq 2.776$.
This gives upper and lower bounds on the determinant of a lattice
$\Gamma $ if one knows $\min(\Gamma )$ and $\min (\Gamma ^*)$:

\begin{lemma} \label{boundsfromgamma} 
Let $\Gamma $ be an $n$-dimensional lattice of minimum $m$ and 
let $r:=\min (\Gamma ^*)$. 
Then for any $b\geq \gamma _n $
$$\left( \frac{m}{b} \right)^n \leq \det(\Gamma ) \leq 
\left( \frac{b}{r} \right) ^n .$$
\end{lemma}

\bew
We have 
$$\frac{m}{\det(\Gamma )^n} \leq b \mbox{ and } r (\det (\Gamma )^n) \leq b $$
since $\det(\Gamma ^*) = \det(\Gamma )^{-1} $.
\eb

\begin{lemma} \label{subsetofdual} 
Let $\Gamma $ be an integral lattice and $X\subset \Gamma ^*$ 
be a linearly independent set with Grammatrix 
$F := ( (x,y) )_{x,y\in X}$. 
Assume that the elementary divisors of $F$ are 
$(\frac{a_1}{b_1},\ldots , \frac{a_k}{b_k} )$ for coprime pairs of 
integers $(a_i,b_i)$. 
Then the product $b_1\ldots b_k $ divides $\det (\Gamma )$.
\end{lemma}

\bew
Let $X = \{ x_1,\ldots , x_k \}$ and $U:=\langle \overline{x}_1,\ldots ,\overline{x}_k \rangle \leq \Gamma ^* / \Gamma $.
Then $U$ is a finite abelian group and any relation matrix 
$A \in \Z^{k\times k}$ such that $\sum _{j=1}^k A_{ij} \overline{x}_j = 0$
for all $i$ satisfies that 
$\sum _{j=1}^k A_{ij} (x_j , x_t)  \in \Z $ for all $i,t$, hence 
$AF \in \Z^{k\times k }$. 
Therefore $\det (A)$ is a multiple of $b_1,\ldots ,b_k$.
\eb

%\subsection{Modular forms.}\label{modforms}

\subsection{Genus symbols and mass formulas.}\label{genus}

We will often encounter the problem to enumerate all even lattices $L$ of 
a given determinant $\det (L) = |L^*/L|$ or of given invariants of 
the abelian group $L^*/L$ (which are the elementary divisors of the 
Grammatrix of $L$).

To list all such genera we implemented a SAGE program
which uses the conditions of Section 7.7 \cite[Chapter 15]{SPLAG} 
in particular \cite[Theorem 15.11]{SPLAG}. 
The python code is available in SAGE and also from 
\cite{Nebehome}.
Also the genus symbol that we use in the classification is the 
one given in \cite[Chapter 15]{SPLAG}.
For odd primes $p$ the lattice $L\otimes \Z _{p}$ has a unique
Jordan decomposition 
$\perp _{i=a}^b p^i f_i $ where the forms $f_i$ are 
regular $p$-adic forms. The determinant $\det (f_i ) = d_i $
is a unit in $\Z_p$ and uniquely determined up to squares.
Let $\epsilon _i := 1$ if $d_i $ is a square and 
$\epsilon _i := -1$ if not. 
Then the $p$-adic symbol of $L$ is 
$((p^i)^{\epsilon_i \dim(f_i) })_{i=a..b} $ where $0$-dimensional 
forms are omitted. 
For $p=2$ we obtain a similar Jordan decomposition which is in 
general not unique. 
To distinguish even  and odd forms $f_i$ we add the oddity, the trace
of a diagonal matrix representing $f_i$, 
as an index, if $f_i$ is odd.
There are also 4 square-classes in $\Z_2^*$, so we put
$\epsilon _i:= 1 $ if $\det(f_i) \equiv \pm 1 \pmod{8} $ 
and $\epsilon _i:= -1 $ if $\det(f_i) \equiv \pm 3 \pmod{8} $.
Since we only use the symbol for even primitive lattices 
the symbol of the first form in the Jordan
decomposition is omitted, as it can be obtained from the others. 
So the symbol for the root lattice $A_5$ is 
$2^{-1}_3 3^1 $.

Any genus $$G(L_1 ) = \{ L \subset (\R^n,(,)) \mid 
L\otimes \Z_p \cong L_1 \otimes \Z_ p \mbox{ for all  primes } p \} $$
of positive definite lattices only consists of finitely many 
isomorphism classes 
$$G(L_1) = [L_1] \cup [L_2] \cup \ldots \cup [L_h ] .$$
The sum 
$$\sum _{i=1}^h  | \Aut (L_i ) | ^{-1 } = \mass (G(L_1)) $$
of the reciprocals of the orders of the automorphism groups of 
representatives $L_i$ of the isomorphism class $[L_i]$ in the 
genus is known as the {\em mass} of the genus an may be computed 
a priori from the genus symbol. 
A formula is for instance given in \cite{mass} and a Magma program 
to calculate its value from the genus symbol
is available from \cite{Nebehome}.
This is used to check whether a list of pairwise non isometric 
lattices in the genus is complete.

To find all lattices in a genus we use the Kneser neighbouring method
\cite{Kneser} (see also \cite{ScharlauHemkemeier}).

\subsection{The maximal even lattices of level dividing 12.}\label{maxeven}

The following table lists all genera of  maximal even lattices $L$ such that 
$\sqrt{12 } L^*$ is again an even lattice. 
The first column gives the genus symbol as explained in \cite[Chapter 15]{SPLAG}, followed by the class number $h$. 
Then we give one representative of the genus which is usually 
a root lattice. $\perp $ denotes the orthogonal sum and 
$(a)$ a 1-dimensional lattice with Grammatrix $(a)$. 
The last column gives the mass of the genus.

$$
\begin{array}{|c|c|c|c|c|}
\hline
\mbox{ genus } & \mbox{ level } & h & \mbox{ repr. } & \mbox{ mass } \\
\hline
3^1 & 3 & 2 & E_6\perp E_8 & 691/( 2^{23} 3^9 5^3 7 \cdot 11\cdot 13)   \\
\hline
2^{2}_6 & 4 & 4 &  D_{14}  & 42151/(2^{25} 3^8 5^3 7^2 11\cdot 13) 
 \\ \hline
2^2_0 3^1 & 12& 8& E_8\perp A_5 \perp (2) &
29713/(2^{17}3^85^37\cdot 11\cdot 13 \\ 
2^{-2} 3^{-1} & 6 &6  & A_2\perp D_{12} &
29713/(2^{24} 3^9 5^2 7 \cdot 11 ) \\ \hline
2^2_2 3^2 & 12 & 28 & A_2 \perp A_2 \perp E_8 \perp (2) \perp (2) &
1683131581/(2^{25} 3^8 5^2 7^2 11 \cdot 13 ) 
 \\
\hline
\end{array} 
$$

\subsection{Designs and strongly perfect lattices}\label{DDD}

For  $m\in \R $, $m>0$ denote by 
$$ S^{n-1}(m) := \{
y\in \R^n \mid (y,y)=m \}$$ 
the $(n-1)$-dimensional sphere of radius $\sqrt{m}$.

\begin{defi}
A finite nonempty set 
$ X \subset S^{n-1}(m)$ is called a {\em spherical $t$-design},
if 
$$\frac{1}{|X|} \sum _{x \in X} f(x) = \int _{S^{n-1}(m)} f(x) d\mu (x) $$
for all polynomials of $f\in \R [x_1,\ldots , x_n ]$ of degree $\leq t$,
where $\mu $ is the  $O(n)$-invariant measure on the sphere, normalized 
such that 
$\int _{S^{n-1}(m)} 1 d\mu (x) = 1 $.
\end{defi}

Since the condition is trivially satisfied for constant polynomials $f$,
and the harmonic polynomials generate the orthogonal complement 
$\langle 1 \rangle ^{\perp } $ with respect to the
$O(n)$-invariant scalar product 
$<f,g>:= \int _{S^{n-1}(m)} f(x) g(x) d \mu (x)  $ on $\R [x_1,\ldots , x_n]$,
it is equivalent to ask that 
$$\sum _{x\in X } f(x) = 0 $$ 
for all harmonic polynomials $f $ of degree $\leq t$.

\begin{defi}{\label{stperf}}
A lattice $\Lambda  \subset \R^n$ is called 
{\em strongly perfect}, if its minimal vectors 
$\Lambda _{\min }  $ form a spherical $4$-design. \\
We call $\Lambda $ {\em dual strongly perfect}, if both lattices,
$\Lambda $ and its {\em dual lattice}
$$\Lambda ^* := \{ v\in \R^n \mid (v,\lambda ) \in \Z \mbox{ for all } \lambda \in \Lambda \} $$ 
are strongly perfect. \\
$\Lambda $ is called {\em universally perfect}, if 
$\theta _{\Lambda , p } = 0 $ for all harmonic polynomials $p$ of degree
2 and 4.
\end{defi}

Let $\Lambda $ be a strongly perfect lattice of dimension $n$,
$m:= \min (\Lambda )$ and choose
$X\subset \Lambda _{m}$ such that $X\cup -X = \Lambda _{m}$
and $X\cap -X = \emptyset$.
Put $s:= |X| = s(\Lambda )$.

By \cite{Venkov} the condition that $\sum _{x\in X} f(x) = 0$ for all
harmonic polynomials $f$ of degree 2 and 4 may be reformulated 
to the condition that for all $\alpha \in \R^n$
$$ (D4) (\alpha ): \ \ \sum _{x\in X} (x,\alpha )^4 = \frac{3sm^2}{n(n+2)} (\alpha , \alpha )^2.$$
Applying the Laplace operator 
$$ \Delta := \sum _{i=1}^n \frac{\partial ^2}{\partial \alpha _i^2 }$$ to $(D4)(\alpha )$  one obtains
$$ (D2) (\alpha ): \ \ \sum _{x\in X} (x,\alpha )^2 = \frac{sm}{n} (\alpha , \alpha ) .$$

Since 
$(x , \alpha )^2 ((x,\alpha )^2-1) $ is divisible by $12$ if $(x,\alpha )\in \Z$
 $$
\frac{1}{12} \sum _{x\in X} (x,\alpha )^2
((x,\alpha )^2-1) = \frac{sm}{12n} (\alpha, \alpha  )
( \frac{3m}{n+2} (\alpha  , \alpha ) - 1 )  \in \Z \ \fa \alpha  \in \Lambda ^* $$

Substituting $\alpha := \xi _1 \alpha _1 + \xi _2 \alpha _2$ 
in  $(D4)$  and comparing coefficients, one finds
$$(D22)(\alpha _1,\alpha _2): \ \ \sum _{x\in X} (x,\alpha _1 )^2 (x,\alpha _2)^2= \frac{sm^2}{n(n+2)} 
(2(\alpha _1, \alpha _2 )^2+(\alpha _1 , \alpha _1)(\alpha _2 , \alpha _2))$$

Also homogeneous polynomials of higher degrees can be
used to obtain linear inequalities. 
To this aim we use the Euclidean inner product on the space of 
all homogeneous polynomials in $n$ variables introduced in \cite{Venkov}.
For $i:=(i_1,\ldots,i_n)\in \Z_{\geq 0}^n$ let 
$x^i := x_1^{i_1} \cdots x_n^{i_n} $ 
and let $c(i):= \frac{(i_1+\ldots + i_n)!}{i_1!\cdots i_n!} $ 
denote the multinomial coefficient.
For two polynomials 
$$f := \sum _{i} c(i)a(i)x^i ,\ 
g := \sum _{i} c(i)b(i)x^i  \mbox{ put } [f,g]:=\sum a(i)b(i) .$$
Let $\omega (x):=(x,x) $ denote the quadratic form.
For $\alpha \in \R^n$ the polynomial $\rho _{\alpha }^m (x):=(x,\alpha )^m $
is homogeneous of degree $m$. 
It is shown in \cite[Proposition 1.1]{Venkov} that 
$$[\rho _{\alpha }^m,f] = f(\alpha ) \mbox{ for all homogeneous } f 
\mbox{ of degree } m $$ 
and by \cite[Proposition 1.2]{Venkov} 
$$ [\omega ^{m/2} , f] = m!\Delta ^{m/2} (f) . $$

\begin{lemma} \label{positiv} 
Let $X\subset \R^n$ be a finite set of vectors of equal norm $m=(x,x)$ such 
that $X\cap -X = \emptyset $.
For $t\in \N $ let $c_t:= (1\cdot 3\cdot 5\cdots (2t-1))/(n(n+2)(n+4) \cdots 
(n+2t-2)) $.
Then 
$$\sum _{x_1,x_2 \in X} (x_1,x_2)^{2t} \geq c_t |X|^2 m^t $$
with equality if and only if $X\disj -X$ is a spherical $2t$-design.
\end{lemma} 

\bew
Define the homogeneous polynomial 
$$p(\alpha ):= \sum _{x\in X}(x,\alpha)^{2t} -c_t |X| (\alpha,\alpha)^t 
= \sum _{x\in X}\rho_x(\alpha)^{2t} -c_t |X| \omega(\alpha)^t .$$
Then the formulas in \cite[Exemple 1.5]{Venkov} show that 
$\Delta ^t(p) = 0$.  
We calculate  the norm 
$$\begin{array}{ll}
 [p,p]  & = \sum _{x\in X} [p,\rho _x^{2t}] - c_t |X| [p,\omega ^t] 
= \sum_{x\in X} p(x) - c_t |X| \Delta ^t(p)   \\
& = \sum _{x_1,x_2\in X} (x_1,x_2)^{2t} -c_t |X|^2 m^t 
\end{array}
$$
which is $\geq 0$ as a norm of a vector in an Euclidean space.
Moreover this sum equals $0$ if $p=0$ which is 
equivalent to the $2t$-design property of the antipodal set
$X\disj -X $.
\eb

\begin{lemma}{\label{linkomb}}(see \cite[Lemma 2.1]{dim10})
Let $\alpha \in \R ^n$ be 
such that $(x , \alpha ) \in \{ 0,\pm 1 ,\pm 2 \} $ for all $x\in X$.
Let $N_2(\alpha ) := \{ x\in X\cup -X  \mid (x,\alpha ) = 2 \}$ and
put $$c:=\frac{sm}{6n}(\frac{3m}{n+2}(\alpha ,\alpha )  -1 ) .$$
Then 
$|N_2(\alpha ) | = c (\alpha , \alpha ) /2 $ and
$$\sum _{x\in N_2(\alpha )} x = c \alpha . $$
\end{lemma}

Lemma \ref{linkomb} will be often applied to $\alpha \in \Lambda ^*$.
Rescale $\Lambda $ such that $\min (\Lambda ) = m = 1$
and let $r:=\min (\Lambda ^*)$.
Since $\gamma(\Lambda ) \gamma (\Lambda ^*) = \min(\Lambda )
\min (\Lambda ^*) \leq \gamma _n ^2$, we 
get $r\leq \gamma _n^2$ and for $\alpha \in \Lambda ^*_r$ we
have $(\alpha , x)^2 \leq r$ for all $x\in \Lambda _1$.
Hence if $r<9$ then $(\alpha , x ) \in \{ 0,\pm 1 ,\pm 2\}$ for all
$x\in X$ and Lemma \ref{linkomb} may be applied.

The next lemma yields  good bounds on $n_2(\alpha )$.

\begin{lemma}{\label{boundn2}}(see \cite[Lemma 2.4]{dim12})
Let $m:= \min (\Lambda ) $ and choose $ \alpha \in \Lambda ^*_r$.
If $r\cdot m < 8$, then
 $$|N_2(\alpha ) |\leq \frac{r m}{8-r m} .$$
\end{lemma}

We also need the case when $r\cdot m = 8$. 
\begin{lemma}\label{boundr8}	
Let $m:= \min (\Lambda ) $ and choose $ \alpha \in \Lambda ^*_r$ such that 
$r\cdot m = 8$.
 $$|N_2(\alpha ) |\leq 2(n-1) .$$
\end{lemma}

\bew
Let $N_2(\alpha ) = \{ x_1,\ldots , x_k \}$ with $k= |N_2(\alpha )| $.
Consider the projections 
$$\overline{x}_i := x_i - \frac{2}{r} \alpha $$ 
of the $x_i $ onto $\alpha ^{\perp }$. 
Then  for all $i\neq j$ 
$$(\overline{x}_i, \overline{x}_j)
 = (x_i,x_j) - \frac{4}{r} 
 = (x_i,x_j) - \frac{m}{2}  \leq 0 $$
since $|(x_i,x_j)  | \leq \frac{m}{2}$ as $x_i,x_j$ are minimal vectors of a 
lattice. 
So the $\overline{x}_i$ form a set of $k$ distinct vectors of
equal length in $(n-1)$-dimensional 
space having pairwise non-positive inner products. 
Moreover $\sum_{i=1}^k \overline{x}_i = 0$ since 
$\sum _{i=1}^k x_i = c \alpha $ by Lemma \ref{linkomb}.
\\
We {\bf claim} that the set 
$$\overline{N}_2(\alpha ) := \{ \overline{x}_1,\ldots , \overline{x}_k \} 
= E_1\disj E_2 \disj \ldots \disj E_{\ell } $$
partitions into $\ell $ disjoint sets such that 
$\sum _{x\in E_i} x = 0 $ for all $i$ and that this is the only
relation. So $\rank (\langle E_i \rangle ) = |E_i| - 1 $ 
and in total 
$$n-1 \geq \rank (\langle \overline{N}_2(\alpha ) \rangle ) = k - \ell 
\geq \frac{k}{2} $$
since $|E_i| \geq 2$ for all $i$. \\
The proof of the claim is standard but for convenience of the reader
we sketch it here.
\\
Any relation between the $\overline{x}_i$ can be written as 
$y:=\sum a_i \overline{x}_i = \sum b_j \overline{x} _j $ 
for nonnegative $a_i, b_j$.
But then 
$$(y,y) = (\sum a_i \overline{x}_i , \sum b_j \overline{x} _j) 
= \sum a_i b_j ( \overline{x}_i,  \overline{x} _j) \leq 0 $$
and hence $y = 0$ so any minimal relation between the
$\overline{x}_i$ has only positive coefficients. 
Subtracting the relations $\sum _{i=1}^k \overline{x}_i = 0$ 
we see that all these coefficients have to be equal. 
\eb

\begin{lemma}{\label{min}}(\cite[Th\'eor\`eme 10.4]{Venkov})
Let $L$  be a strongly perfect lattice of dimension $n$.
Then 
$$\gamma(L) \gamma (L^*) =  \min(L) \min(L^*) \geq \frac{n+2}{3} .$$
A strongly perfect lattice $L$ is called 
{\em of minimal type}, if $\min(L) \min(L^*) = \frac{n+2}{3}$ and 
{\em of general type} otherwise. 
\end{lemma}

\subsection{Certain sublattices.}

The next lemma about indices of sublattices are used quite often
in the argumentation below.
Since we are dealing with norms modulo some prime number $p$, we may
pass to the localization
$\Z _p := ( \Z \setminus p\Z )^{-1} \Z \subset \Q $ of
$\Z $ at $p$.

\begin{lemma}\label{sublat}(see \cite[Lemma 2.8, 2.9]{dim12})
\\
{\bf a)} 
Let $\Gamma $ be a $\Z _2$-lattice such that $(\gamma ,\gamma ) \in \Z _2 $ for 
all $\gamma \in \Gamma $.
Let $\Gamma ^{(e)} := \{\alpha \in \Gamma \mid (\alpha,\alpha ) \in 2 \Z _2 \}$.
If $\Gamma ^{(e)} $ is a sublattice of $\Gamma $, then 
$[\Gamma : \Gamma ^{(e)}] \in \{ 1,2,4 \}$.
\\
{\bf b)}
Let $\Gamma $ be a $\Z _3$-lattice such that $(\gamma ,\gamma ) \in \Z _3 $ for 
all $\gamma \in \Gamma $.
Let $\Gamma ^{(t)} := \{ \alpha \in \Gamma \mid (\alpha,\alpha ) \in 3 \Z _3 \}$.
Assume that 
$$ (\alpha ,\beta )^2 - (\alpha , \alpha ) (\beta , \beta)  \in 3\Z_3 \mbox{ for
all } \alpha , \beta \in \Gamma .$$
Then 
$\Gamma ^{(t)} $ is a sublattice of $\Gamma $ and
$[\Gamma : \Gamma ^{(t)}] \in \{ 1,3 \}$.
\end{lemma}

Note that in both cases of Lemma \ref{sublat}, if $\Gamma $ is 
universally perfect, then so is $\Gamma ^{(e)}$ resp. $\Gamma ^{(t)}$ 
since both lattices consist of certain layers of $\Gamma $.

\section{General type.}

\subsection{Kissing numbers.}

Let $\Lambda $ be a strongly perfect lattice in dimension 14, rescaled
such that $\min (\Lambda ) = 1$.
Then by Lemma \ref{min}  and the Cohn-Elkies bound that 
$\gamma _{14} \leq 2.776 $ in \cite{Elkies} we find 
$$ 16/3 \leq r:=\min(\Lambda ^* ) \leq \gamma _{14}^2 \leq 7.71 .$$
Hence $\alpha \in \Lambda ^* _r$ satisfies the hypothesis of  Lemma \ref{linkomb}
and Lemma \ref{boundn2} yields 
$$n_2(\alpha ) \leq \frac{r}{8-r} \leq \frac{\gamma _{14}^2}{8-\gamma _{14}^2}
\leq 26.3 < 27 .$$

Let $s:= |X|$ where $X\disj -X = \Lambda _1$ be half the 
kissing number of $\Lambda $.
Then by the bound given in \cite{Anstreicher}
$\frac{15\cdot 14}{2} = 105 \leq s \leq 1746 $ and 
$r < 7.71 $ is a rational solution of 
$$ 1 \neq n_2(\alpha ) = \frac{sr}{12\cdot 14\cdot 16} (3r-16)  \leq \frac{r}{8-r}$$
and $\frac{s\cdot r}{14} \in \Z$, $\frac{3sr^2}{14\cdot 16} \in \Z $.
Going through all possibilities by a computer we find:

\begin{prop}{\label{poss14}}
With the notation above, one of the following holds:
$$
\begin{array}{|c|c|c|c|c|c|c|c|c|c|c|}
\hline
n_2 & 2 & 2 & 2 &  3 & 3 &  4 &4  & 4 & 4  \\
\hline
s &  324 & 448 & 1200 & 486 & 672   & 225 & 343 & 363 & 525 \\
\hline
r & 56/9 & 6 & 28/5 & 56/9 & 6  &112/15& 48/7 & 224/33 & 32/5 \\
\hline
\end{array}
$$
$$
\begin{array}{|c|c|c|c|c|c|c|c|c|c|}
\hline
n_2 & 5 & 5 & 8 &  8 & 11 &  12 & 19  & 20   \\
\hline
s &  384 & 504 & 450 & 567 & 672   & 675 & 968 & 1029  \\
\hline
r & 7 & 20/3 & 112/15 & 64/9 & 22/3  &112/15& 84/11 & 160/21  \\
\hline
\end{array}
$$
or $\Lambda $ is of minimal type, i.e. $r=16/3$.
\end{prop}

Most of the cases in Proposition \ref{poss14} are ruled out quite easily.
More precisely we will prove the following
\begin{theorem}{\label{theorem1}}
Let $\Lambda $ be a strongly perfect lattice of dimension $14$, 
with $\min (\Lambda ) := 1 $.
Let $s:= \frac{1}{2} | \Lambda _1 |$  and $r:= \min (\Lambda ^* )$.
Then we have the following three possibilities
\begin{itemize}
\item[(a)] $s=672$ and $r=6$ or $r=22/3$.
\item[(b)] $s=504$ and $r=20/3$.
\item[(c)] $r=16/3$, which means that $\Lambda $ is of minimal type.
\end{itemize}
\end{theorem}

\subsection{Proof of Theorem \ref{theorem1} for $s\neq 450$.}

For the proof we scale $\Lambda $, such that $\min (\Lambda ) = 1$.
For $\alpha \in \Lambda ^*$  write
$(\alpha , \alpha ) = \frac{p}{q}$ with coprime integers $p$ and $q$.
Then
$$(\star ) \ \ \frac{1}{12} (D4-D2)(\alpha ) = \frac{s}{2^73\cdot 7} \frac{p}{q^2}
(3p-16q)  \in \Z .$$
Moreover $$
\begin{array}{ll}
D2(\alpha ) & = \frac{s}{14} \frac{p}{q}  \in \Z  \\
D4(\alpha ) & = \frac{3s}{14\cdot 16} \frac{p^2}{q^2}  \in \Z .
\end{array} $$

\begin{lemma}
$(s,r) \neq (324, 56/9) $. 
\end{lemma}

\bew
Assume that $s=324=2^23^4$ and $r=56/9$.
Then $D4$ yields that $\frac{3^5}{2^37} \frac{p^2}{q^2} \in \Z $ hence 
$q$ divides $3^2$ and $7$ divides $p$. 
Moreover by $(\star )$ 
$\frac{3^3}{2^57} \frac{p}{q^2} (3p-16q) \in \Z $ yields that $p$ is 
divisible by $2^3$.
Therefore 
$\Gamma := \sqrt{\frac{9}{28}} \Lambda ^*$ is an even lattice of minimum 2 and
$\Gamma ^*$ has minimum $\frac{28}{9} > 2$ which is a contradiction.
\eb

\begin{lemma}
$(s,r) \neq (448, 6) $. 
\end{lemma}

\bew
Assume that $s=448=2^67$ and $r=6$.
Then $D4$ yields that $6 \frac{p^2}{q^2} \in \Z $ hence 
all norms in $\Lambda ^*$ are integral ($q=1$).
Moreover by $(\star )$ 
$\frac{1}{6} p (3p-16) \in \Z $ yields that $p$ is 
divisible by $6$.
Therefore 
$\Gamma := \sqrt{\frac{1}{3}} \Lambda ^*$ is an even lattice of minimum 2 and
$\Gamma ^*$ has minimum $3 > 2$ which is a contradiction.
\eb

\begin{lemma}
$(s,r) \neq (1200, 28/5) $. 
\end{lemma}

\bew
Assume that $s=1200=2^43\cdot 5^2$ and $r=28/5$.
Then $D4$ yields that $\frac{3^25^2}{14} \frac{p^2}{q^2} \in \Z $ hence 
$q$ divides 15 and 14 divides $p$.
Moreover by $(\star )$ 
$\frac{5^2}{2^37} \frac{p}{q^2} (3p-16q) \in \Z $ shows that 
 $p$ is a multiple of $4$.
Put
$\Gamma := \sqrt{\frac{5}{14}} \Lambda ^*$.
Then $\min (\Gamma ^*) = \frac{14}{5} $ and equation (D22) shows that 
$$3\cdot 14 (2(\alpha,\beta )^2 + (\alpha,\alpha)(\beta,\beta ) ) \in \Z
\mbox{ for all } \alpha,\beta \in \Gamma $$
Therefore 
$$\Gamma ^{(t)}:=\{ \alpha \in \Gamma \mid (\alpha,\alpha ) \in \Z \} =
\Gamma \cap \Gamma ^* $$
is an even sublattice of $\Gamma $ of index $1$ or $3$.
Moreover $\min(\Gamma ^{(t)}) = 2$ and $\min (\Lambda ) = \frac{14}{5} > 2$
contradicting the fact that $\Gamma ^{(t)} \subset \Lambda $.
\eb

\begin{lemma}
$(s,r) \neq (486, 56/9) $. 
\end{lemma}

\bew
Assume that $s=486=2\cdot 3^5$ and $r=56/9$.
Then $D4$ yields that $\frac{3^6}{2^47} \frac{p^2}{q^2} \in \Z $ hence 
$q$ divides $3^3$ and 28 divides $p$.
Moreover by $(\star )$ 
$\frac{3^4}{2^67} \frac{p}{q^2} (3p-16q) \in \Z $ yields that $q$ divides $9$  
and $p$ is
divisible by $2^3$.
Therefore 
$\Gamma := \sqrt{\frac{9}{28}} \Lambda ^*$ is an even lattice of minimum 2 and
$\Gamma ^*$ has minimum $\frac{28}{9} > 2$ which is a contradiction.
\eb

\begin{lemma}
$(s,r) \neq (225, 112/15) $. 
\end{lemma}

\bew
Assume that $s=225= 3^25^2$ and $r=112/15$.
Then $D4$ yields that $\frac{3^35^2}{2^57} \frac{p^2}{q^2} \in \Z $ hence 
$q$ divides $15$ and 56 divides $p$.
Moreover by $(\star )$ 
$\frac{3\cdot 5^2}{2^77} \frac{p}{q^2} (3p-16q) \in \Z $ yields that
$p$ is
divisible by $2^4$.
Therefore 
$\Gamma := \sqrt{\frac{15}{56}} \Lambda ^*$ is an even lattice of minimum 2 and
$\Gamma ^*$ has minimum $\frac{56}{15} > 2$ which is a contradiction.
\eb

\begin{lemma}
$(s,r) \neq (343, 48/7) $. 
\end{lemma}

\bew
Assume that $s=343= 7^3$ and $r=48/7$.
Then $D4$ yields that $\frac{3\cdot 7^2}{2^5} \frac{p^2}{q^2} \in \Z $ hence 
$q$ divides $7$ and 8 divides $p$.
Moreover by $(\star )$ 
$\frac{7^2}{2^73} \frac{p}{q^2} (3p-16q) \in \Z $ yields that
$p$ is
divisible by $2^43$.
Therefore 
$\Gamma := \sqrt{\frac{7}{24}} \Lambda ^*$ is an even lattice of minimum 2 and
$\Gamma ^*$ has minimum $\frac{24}{7} > 2$ which is a contradiction.
\eb

\begin{lemma}
$(s,r) \neq (363, 224/33) $. 
\end{lemma}

\bew
Assume that $s=363= 3\cdot 11^2$ and $r=224/33$.
Then $D4$ yields that $\frac{3^2\cdot 11^2}{2^57} \frac{p^2}{q^2} \in \Z $ hence 
$q$ divides $33$ and $2^37$ divides $p$.
Moreover by $(\star )$ 
$\frac{11^2}{2^77} \frac{p}{q^2} (3p-16q) \in \Z $ yields that
$p$ is
divisible by $2^4$.
Therefore 
$\Gamma := \sqrt{\frac{33}{56}} \Lambda ^*$ is an even lattice of minimum 4 and
$\Gamma ^*$ has minimum $\frac{56}{33} $.
In the new scaling the equation 
$D22$ reads as 
$\frac{14}{3} (2 (\alpha , \beta )^2 + (\alpha ,\alpha) (\beta ,\beta ) ) \in \Z $ 
for all $\alpha, \beta \in \Gamma $.
Therefore
$$\Gamma ^{(t)}:= \{ \gamma \in \Gamma \mid  (\gamma,\gamma ) \in 3 \Z \}$$
is a sublattice of index 3 in $\Gamma $ and hence $3^{12}$ divides 
$\det (\Gamma )$. 
This yields 
$$ 531441 = 3^{12} \leq \det(\Gamma ) \leq (\frac{\gamma _{14}}{\min (\Gamma ^*)})^{14} 
 < 982 $$
which is a contradiction.
\eb

\begin{lemma}
$(s,r) \neq (525, 32/5) $. 
\end{lemma}

\bew
Assume that $s=525= 3\cdot 5^27$ and $r=32/5$.
Then $D4$ yields that $\frac{3^2\cdot 5^2}{2^5} \frac{p^2}{q^2} \in \Z $ hence 
$q$ divides $15$ and $2^3$ divides $p$.
Moreover by $(\star )$ 
$\frac{5^2}{2^7} \frac{p}{q^2} (3p-16q) \in \Z $ shows that 
$p$ is a multiple of $2^4$.
Put
$\Gamma := \sqrt{\frac{5}{8}} \Lambda ^*$.
Then $\min (\Gamma ) = 4$ and 
$\Gamma ^* = \frac{\frac{8}{5}} \Lambda $ has minimum $\frac{8}{5} $.
The equality (D22) yields that 
$$12((\alpha,\beta)^2 + (\alpha,\alpha )(\beta,\beta) ) \in \Z 
\mbox{ for all } \alpha,\beta \in \Gamma $$
hence 
$$\Gamma ^{(t)}:=\{\alpha \in \Gamma \mid (\alpha,\alpha ) \in \Z \}
=\Gamma \cap \Gamma ^* $$ is a sublattice of $\Gamma $ of index 3 or 1.

Let $\alpha \in \Min(\Gamma )$. Then $(\alpha ,\alpha ) = 4$ and hence
$\alpha \in \Gamma ^{(t)}$.
Moreover
$$|N_2(\alpha ) | = \frac{1}{12} (D4-D2)(\alpha ) = \frac{1}{2} (\alpha ,\alpha ) 
(3 (\alpha , \alpha ) - 10) = 4 $$
and  $N_2(\alpha )  =: \{ x_1,x_2 := \alpha - x_1 ,x_3,x_4:= \alpha - x_3 \} $ 
for certain $x_1, x_3 \in \Min (\Gamma ^*)$.
Since these are minimal vectors of a lattice, we have 
$(x_1,x_3) \leq \frac{1}{2} (x_1,x_1) = \frac{4}{5} $ and hence
$(x_1,x_4) = (x_1,\alpha - x_3) = 2 - (x_1,x_3) \geq 2 - \frac{4}{5} = \frac{6}{5} $
yielding the vector $x_1-x_4 \in \Gamma ^*$ of norm $<\min (\Gamma ^*)$, a 
contradiction.
\eb

\begin{lemma}
$(s,r) \neq (384, 7) $. 
\end{lemma}

\bew
Assume that $s=384= 2^73$ and $r=7$.
Then $D4$ yields that $\frac{2\cdot 3^2}{7} \frac{p^2}{q^2} \in \Z $ hence 
$q$ divides $3$ and $7$ divides $p$ and therefore
$\Gamma := \sqrt{\frac{6}{7}} \Lambda ^*$ is an even lattice of minimum 6 and
$\Gamma ^*$ has minimum $\frac{7}{6} $.
In the new scaling the equation 
$D22$ reads as 
$\frac{7}{3} (2 (\alpha , \beta )^2 + (\alpha ,\alpha) (\beta ,\beta ) ) \in \Z $ 
for all $\alpha, \beta \in \Gamma $.
Therefore
$$\Gamma ^{(t)}:= \{ \gamma \in \Gamma \mid  (\gamma,\gamma ) \in 3 \Z \}$$
is a sublattice of index 3 in $\Gamma $ and hence $3^{12}$ divides 
$\det (\Gamma )$. 
This yields 
$$ 531441 = 3^{12} \leq \det(\Gamma ) \leq (\frac{\gamma _{14}}{\min (\Gamma ^*)})^{14} 
 < 186474 $$
which is a contradiction.
\eb

\begin{lemma}
$(s,r) \neq (567, 64/9) $. 
\end{lemma}

\bew
Assume that $s=567= 3^47$ and $r=64/9$.
Then $D4$ yields that $\frac{ 3^5}{2^5} \frac{p^2}{q^2} \in \Z $ hence 
$q$ divides $3^2$ and $2^3$ divides $p$.
Moreover by $(\star )$ 
$\frac{3^3}{2^7} \frac{p}{q^2} (3p-16q) \in \Z $ yields that
$p$ is
divisible by $2^4$.
Therefore 
$\Gamma := \sqrt{\frac{9}{8}} \Lambda ^*$ is an even lattice of minimum 8 and
$\Gamma ^*$ has minimum $\frac{8}{9} $.
Choose $\alpha \in \Min (\Gamma )$.
Then $N_2(\alpha ) = \{ x_1,\ldots , x_8 \}$ with $\sum _{i=1}^8 x_i = 2 \alpha $.
Moreover 
$$ 4 = (x_1 , 2\alpha ) = ( x_1 , \sum _{i=1}^8 x_i ) = (x_1,x_1) + \sum _{i=2}^8 (x_1,x_i) 
\leq \frac{8}{9} + 7 \frac{4}{9} = 4 $$
yields that 
$N_2(\alpha ) = \frac{2}{3}  A_8 $. 
In particular the denominator of the determinant of the Grammatrix of 
$N_2(\alpha )$ is $9^7$. 
By Lemma \ref{subsetofdual}
this implies that 
$9^7$ divides the order $|\Gamma ^* / \Gamma | = \det (\Gamma ) \leq 
(\gamma _{14}/ (8/9)) ^{14} < 2 \cdot 9^7  $.
Hence $\det (\Gamma ) = 9^7$ and 
$$ L := \langle \Gamma , 3 x_1,\ldots , 3 x_7 \rangle $$ is an even
overlattice containing $\Gamma $ of index $3^7$.
Therefore $L$ is unimodular, a contradiction since there is no 
even unimodular lattice in dimension 14.
\eb

\begin{lemma}
$(s,r) \neq (675,112/15)$.
\end{lemma}

\bew
Assume that $s=675= 3^35^2$ and $r=112/15$.
Then $D4$ yields that $\frac{ 3^45^2}{2^57} \frac{p^2}{q^2} \in \Z $ hence 
$q$ divides $3^25$ and $2^37$ divides $p$.
Moreover by $(\star )$ 
$\frac{3^25^2}{2^77} \frac{p}{q^2} (3p-16q) \in \Z $ yields that $q$ 
divides  $15$ and $p$ is divisible by $2^4$.
Therefore 
$\Gamma := \sqrt{\frac{15}{56}} \Lambda ^*$ is an even lattice of minimum 2 and
$\Gamma ^*$ has minimum $\frac{56}{15} > 2 $
which is a contradiction.
\eb

\begin{lemma}
$(s,r) \neq (968,84/11)$.
\end{lemma}

\bew
Assume that $s=968=2^311^2$ and $r=84/11$.
Then $D4$ yields that $\frac{ 3\cdot 11^2}{2^27} \frac{p^2}{q^2} \in \Z $ hence 
$q$ divides $11$ and $2\cdot 7$ divides $p$.
Moreover by $(\star )$ 
$\frac{11^2}{2^43\cdot 7} \frac{p}{q^2} (3p-16q) \in \Z $ yields that 
 $p$ is divisible by $2^23\cdot 7$.
Therefore 
$\Gamma := \sqrt{\frac{11}{42}} \Lambda ^*$ is an even lattice of minimum 2 and
$\Gamma ^*$ has minimum $\frac{42}{11} > 2 $
which is a contradiction.
\eb

\begin{lemma}
$(s,r) \neq (1029,160/21)$.
\end{lemma}

\bew
Assume that 
$(s,r) = (1029,160/21)$ and choose $\alpha \in \Min (\Lambda ^*)$.
Then $N_2(\alpha )$ generates a rescaled version of the root lattice $A_{20}$,
which is of course impossible in dimension 14.
\eb
\subsection{The case $s=450$.}

To show Theorem \ref{theorem1} it suffices to treat the case $s=450$ 
which is done in this subsection because it involves a bit more 
calculations than the cases treated in the previous section.

Let $\Lambda $ be a strongly perfect lattice of dimension 14 and 
assume that $s(\Lambda ) = \frac{1}{2} |\Lambda _{min} | = 450 $.
Without loss of generality we assume that $\Lambda $ is generated by
its minimal vectors.
Rescale such that $\min (\Lambda ) = \frac{28}{15}$ and let 
$\Gamma := \Lambda ^*$. 
Then $\min (\Gamma ) = 4$ and 
the equations $(D4)$ and $\frac{1}{12} ((D4)-(D2)) $ yield that 
$$ 21 (\alpha , \alpha )^2 \in \Z  \mbox{ and } 
 \frac{1}{4} (\alpha , \alpha ) (7 (\alpha , \alpha ) - 20) \in \Z  
$$
hence $\Gamma $ is an even lattice.
Since $\Lambda $ is generated by elements of norm $\frac{28}{15}$, 
we have that 15 divides $\det (\Gamma )$ and hence by Lemma \ref{boundsfromgamma}
$$ (+) \ \ \ 11 \cdot 15 \leq \det (\Gamma ) \leq 17 \cdot 15 .$$
Moreover for any prime  $p\neq 3,5$ that divides $\det (\Gamma )$ 
the Sylow $p$-subgroup of $\Gamma ^*/\Gamma $ is generated by 
isotropic elements and hence cannot be cyclic. 
In particular
$p^2$ divides $\det (\Gamma )$ and therefore 
$\det (\Gamma ) = 2^a 3^b 5^c $ for some $a,b,c \in \N_0$.
Let $M\supseteq \Gamma $ be a maximal even overlattice of
$\Gamma $. Then also $\min (M^*) \geq \frac{28}{15}$. 
We have the following possibilities for the genus of $M$:

$$
\begin{array}{|r|r|r|}
\hline 
 \mbox{ genus } & h  & \mbox{ mass } \\ 
\hline 
 3^1 &      2  & % E6+E8
 691 / 
 ( 2^{23} 3^9 5^3 7\cdot 11 \cdot 13  ) \sim 
 3.344666\cdot 10^{-14}
 \\ 
 \hline
 3^{-1} 5^1 & 8 & %b E8+A4+A2
 650231/
 ( 2^{21} 3^8 5^3 7^2 13  ) \sim 
 5.9349666\cdot 10^{-10}
\\
3^1 5^{-1} &  9 &  %a A14
 650231/
 ( 2^{21} 3^8 5^3 7^2 13  ) \sim 
 5.9349666\cdot 10^{-10}
  \\
 \hline
3^{-1} 5^{-2} & 48 & %A2 + Qinfty,5 (normform) + E8
5407504111/
 ( 2^{23} 3^9 5^3 7 \cdot 11  ) \sim 
3.4026316\cdot 10^{-6}
\\
\hline
%2^2_6 & &
 %\\
%\hline
%2^{-2}_4 3^1 & & % A5+(2) + E8 aber nicht anistrop
%\\
%2^{-2}_0 3^{-1} & &  % A2+D12
%\\
%\hline
%2^2_2 5^1 & & 
%\\
%2^{-2} _2 5^{-1} & & 
%\\
\hline
2^2_0 3^1 5^{-1} & 93 &  %A9+A5, a
82579337/
 ( 2^{15} 3^8 5^3 7^2 13  ) \sim 
4.823941\cdot 10^{-6}
\\
2^2_0 3^{-1} 5^{1} & 91 &  % E8 + A4 + (2,6) , c
82579337/
 ( 2^{15} 3^8 5^3 7^2 13  ) \sim 
4.823941\cdot 10^{-6}
\\
2^{-2}_0 3^1 5^1 & 46 &  %D4+A4+E6, b, dasselbe wie 2A_2+A_4+E_8
82579337/
 ( 2^{22} 3^7 5^3 7 \cdot 13  ) \sim 
7.914278\cdot 10^{-7}
\\
2^{-2}_0 3^{-1} 5^{-1} & 48 &  %d Obergitter von sqrt{2} A4 + A2 + E8
82579337/
 ( 2^{22} 3^7 5^3 7 \cdot 13  ) \sim 
7.914278\cdot 10^{-7}
\\
%\\
%\hline
%2^{-2}_2 3^2 5^1 & & 
%\\
\hline
 \end{array} 
 $$

The table gives the genus symbol as explained in Section \ref{genus}, the
class number $h$ of the genus followed by the mass.

\bew
We show that the list of possible elementary divisors of $M$ is complete,
then the list of possible genus symbols is obtained with the SAGE 
program mentioned in Section \ref{genus}. 

\knubbel Assume first that the determinant of $M$ is $3^b5^c$. 
Then by Milgram's formula (see  \cite[Cor. 5.8.2 and 5.8.3]{Scharlau}, \cite[Lemma 2.3, 2.4]{dim10}) $b$ is odd and hence 
$\det(M) = 3$, $3\cdot 5$, or $3\cdot 5^2$. 

\knubbel Assume now that $\det(M)$ is even. Then the Sylow 2-subgroup of 
$M^*/M$ has rank 2 and is anisotropic.
Since the Sylow 2-subgroup of $\Gamma ^*/\Gamma $ is generated by 
isotropic elements, the index of $\Gamma $ in $M$ is even.
By the bound in $(+)$, we obtain $\det(\Gamma ) = 2^4 15 $ and hence
$\det(M) = 2^2 \cdot 3 \cdot 5 $. 
\eb

For none of the 345 relevant maximal even lattices $M$ the dual lattice
has minimum $\geq \frac{28}{15}$ so there is no strongly perfect lattice
$\Lambda $
of dimension 14 with kissing number $2\cdot 450$ and 
$\min(\Lambda ) \min (\Lambda ^*) = 112/15 $. 
\section{Dual strongly perfect lattices of general type.}

In this section we show the following theorem.

\begin{theorem}
Let $\Lambda $ be a 14-dimensional dual strongly perfect lattice.
Then $\Lambda $ is of minimal type.
\end{theorem}

To prove this theorem it is enough to consider the 
three remaining cases
($s=672$ and $\gamma =6$ or $\gamma =22/3$ resp. 
$s=504$ and $\gamma =20/3$)
 of Theorem \ref{theorem1}.

So let $\Lambda $ be a dual strongly perfect lattice that is not of 
minimal Type.

\begin{lemma}
$(s,\gamma ) \neq (672,22/3) $.
\end{lemma}

\bew
Let $\Gamma := \Lambda ^*$ and scale such that 
$$\min (\Gamma ) =22, \min (\Lambda ) = \frac{1}{3}. $$
Then for all $\alpha \in \Gamma $
$$\frac{1}{12} (D4-D2) = \frac{1}{12} (\alpha,\alpha ) ( (\alpha,\alpha) -16 ) 
\in \Z $$
implies that $\Gamma $ is an even lattice in which all norms are 
$\equiv 1, 0 \pmod{3} $. 
The equation $(D22)$ implies that for all $\alpha , \beta  \in \Gamma $ 
$$\frac{1}{3} (2(\alpha , \beta )^2 + (\alpha,\alpha ) (\beta , \beta ) ) \in \Z $$ and hence by Lemma \ref{sublat} the lattice
$$\Gamma ^{(t)} = \{ \alpha \in \Gamma \mid 3 \mid (\alpha,\alpha ) \} $$
is a sublattice of $\Gamma $ of index 3.
In particular $3^{12}$ divides the determinant of $\Gamma $ so write
$\det (\Gamma ) = 3^{12} d$.

Now $\Lambda $ is dual strongly perfect, so also $\Gamma $ is 
a strongly perfect lattice of dimension 14. By Theorem \ref{theorem1} 
the lattice $\Gamma $ has the same parameters as $\Lambda $.
In particular $L:=\sqrt{66} \Lambda $ is again an
even lattices of minimum 22 possessing a sublattice $L^{(t)}$ of
index 3 in which all inner products are multiples of 3.
Hence $3^{12} $ divides $$\det (L) = \frac{3^{14}}{22^{14}} \det (\Lambda ) 
= 3^2\cdot 22^{14} \cdot \frac{1}{d} $$ 
since $\Lambda = \Gamma ^*$. 
But $d$ is an integer and therefore we obtain a contradiction.
\eb

\begin{lemma}\label{beispmodular}
$(s,\gamma ) \neq (504,20/3) $.
\end{lemma}

\bew
Let $\Gamma := \Lambda ^*$ and scale such that 
$$\min (\Gamma ) =10, \min (\Lambda ) = \frac{2}{3}. $$
Then $\frac{1}{12} (D4-D2)$ shows that $\Gamma $ is an even lattice.
For $\alpha \in \Gamma _{10}$ consider 
$$N_2(\alpha ) = \{ x\in \Lambda \mid (x,x) = \frac{2}{3} , (x,\alpha ) = 2 \} .$$
Then $|N_2(\alpha ) | = 5$ and $\sum _{x\in N_2(\alpha )} x = \alpha $ 
which implies that $(x,y) = \frac{1}{3}$ for all $x\neq  y\in N_2(\alpha)$ and
hence 
$\langle N_2(\alpha ) \rangle \cong \frac{1}{\sqrt{3}} A_5$.
In particular $3^4$ divides $\det(\Gamma )$.

Moreover by Theorem \ref{theorem1}, the lattice $\Gamma $ has the
same parameters as $\Lambda $ and therefore 
$\sqrt{15} \Lambda $ is again even of minimum 10 and hence 
$\det(\Gamma) = 3^a 5^b $ for some $a\geq 4$.  
Also by Milgram's formula (see \cite[Cor. 5.8.2 and 5.8.3]{Scharlau}) 
we know that the exponent $a$ of $3$ is odd. 
Interchanging the role of $\Lambda $ and $\Gamma $ if necessary we may 
assume that $b\geq 7$. 
By Lemma \ref{boundsfromgamma}
$$61962301 \leq \det(\Gamma ) \leq 471140124 $$ 
so $(a,b)$ is one of $(7,7)$, $(5,8)$. 
For both pairs $(a,b)$ we have 2 possible genera of even lattices. 
$$ \begin{array}{|c|c|}
\hline
\mbox{ genus } & \mbox{ mass }  \\
\hline
3^{7} 5^7 &  > 1 043 713 033 837 \\
\hline
3^{-7} 5^{-7} &  > 1 043 713 033 837 \\
\hline
3^{5} 5^{8} &  > 52297468940 \\
\hline
3^{-5} 5^{-8} &  > 52130384375 \\
\hline
\end{array}
$$

Note that the first two genera are dual to each other.
In view of the masses of the genera given in the table, it is infeasible 
to list all lattices in these genera. 
Instead we will construct one lattice $L$ in each of these genera as a sublattice
of a suitable maximal even lattice and then use explicit calculations 
in the spaces of modular forms 
generated by the theta series of $L$ and the cusp forms for $\Gamma _0(12)$ 
to the relevant character $\left( \frac{-\det(L)}{\cdot } \right) $:

Representatives of the genera of maximal even lattices of 
exponent dividing 15 are 
$A_{14} $ ($3^15^{-1}$), $A_2\perp A_4 \perp E_8$ ($3^{-1} 5^1$), 
$E_6 \perp E_8 $ ($3^1$), $A_2 \perp E_8 \perp L_4$ ($3^{-1} 5^{-2} $) 
where $L_4$ is the 4-dimensional maximal even lattice of determinant $5^2$ 
with Grammatrix 
$F = \left( \begin{array}{cccc} 
2&1&1&1 \\
1&2&0&1 \\
1&0&4&2 \\
1&1&2&4 \end{array} \right) $.
One constructs representatives $L_1,\ldots , L_4$ 
of the 4 relevant genera of lattices 
of determinant $3^75^7$ respectively $3^55^8$ as sublattices of these 
4 lattices above.

We now use the explicit knowledge of the spaces of modular forms for 
$\Gamma _{0}(15)$ with the character given by the Kronecker symbol of 
the discriminant of $L_i$. 
Since the strategy is the same in all cases, we explain it in the 
case where $\Gamma $ is a sublattice of index $3^35^3$ of $A_{14}$.
The relevant space of modular forms 
${\cal M}_{7} (\Gamma _{0}(15), \left( \frac{-15}{\cdot } \right) ) $
has dimension 14 and its cuspidal subspace ${\cal S}$ has dimension 10.
We construct the theta series $T$ of some lattice in the genus of 
$\Gamma $. Then we know that $\theta _{\Gamma } = T + f $ for some
$f \in {\cal S}$.
Moreover the $q$-expansion of $\theta _{\Gamma }$ starts with 
$$1 +  0 q^2 + 0 q^4 + 0 q^6 + 0 q^8 + 2\cdot 21 \cdot 24 q^{10}  + \ldots .$$
Applying the Atkin-Lehner involution $W_{15}$ we obtain the 
same conditions on $W_{15} (\theta _{\Gamma })$. 
To do this in Magma, we have to multiply $\theta _{\Gamma }$ with the
theta series $f_1$ of the lattice with Grammatrix $\left( \begin{array}{cc}
4 & 1 \\ 1 & 4 \end{array} \right) $ to obtain a modular form for 
$\Gamma _0(15)$ of Haupttypus (without a character), that is then of
weight 8. To such forms, we may apply the Atkin-Lehner operator in 
Magma. Then $$W_{15}(\theta_{\Gamma }) = W_{15}(\theta_{\Gamma} f_1 )/ f_1 $$
(since all lattices of dimension 2 are similar to their dual lattices).
This way we obtain 12 linear conditions on the cusp form $f \in {\cal S}$, 
of which only 9 are linearly independent. 
So there is a 1-parametric space of solutions 
$\{ \theta + a f_0 \mid a\in \Z \} $ 
where 
$$\begin{array}{ccrrrrrrr}
\theta = & 1  & +1008q^{10}  & +1896q^{12}  & +43124q^{14}  & -210044q^{16}  & +340244q^{18}  & +755692q^{20} & -  \ldots \\
f_0 = & & &  q^{12} & - 11q^{14} & + 44q^{16} & - 51q^{18} & - 154q^{20} &  +  \ldots 
\end{array} 
$$
One immediately sees that this space does not contain a function with 
positive coefficients, since 
$210044/44 = 52511/11  > 43124/11 .$
So no lattice in the genus of $L$ has minimum 10, kissing number 1008 
such that its rescaled dual lattice $\sqrt{15} L^*$ also has
minimum 10 and kissing number 1008. 

Similarly the other three genera are treated, where one gets a unique solution
(with negative coefficient at $q^{14}$) in the two cases where the 
determinant of $\Gamma $ is $3^5 5^8 $. 
\eb

\begin{lemma} 
If $\Lambda $ is a dual strongly perfect lattice of dimension 14 then
$(s(\Lambda ), \gamma '(\Lambda )^2 ) \neq (672, 6) $.
\end{lemma} 

\bew
Let $\Lambda $ be such a lattice and let $\Gamma := \Lambda ^*$ rescaled
such that $\min (\Gamma ) = 6$, $\min (\Lambda ) = 1$.
Then the norms of the elements in $\Gamma $ are in $\frac{2}{3} \Z $ and 
the equation $D22$ shows that $\Gamma $ contains a sublattice
$$\Gamma ^{(t)} := \{ \alpha \in \Gamma \mid 
(\alpha,\alpha ) \in \Z \} = \Lambda \cap \Gamma $$ 
of index 1 or 3. 
Since also $\Gamma $ is strongly perfect, we may rescale to see that 
an analogous property is hold by $\sqrt{6}\Lambda $.
Let $Z\disj -Z:= \Min (\Gamma )$ and $X\disj -X = \sqrt{6} \Min (\Lambda )$.
Then $T:=\pm X\cup Z$ is a spherical 
5-design consisting of $2\cdot 2 \cdot 672 $ vectors of norm 6.
Moreover the inner products in $T$ are 
$$ 0,\pm 1, \pm 2, \pm 3, \pm 6, \pm \sqrt{6}, \pm 2 \sqrt{6} .$$
We now choose some $\alpha \in Z$  and choose signs of the elements in 
$Z$ such that $(z,\alpha ) \geq 0 $ for all $z\in Z$. Put
$$M_i (\alpha ) := \{ z\in Z \mid (z,\alpha ) = i \} \mbox{ and } 
m_i := |M_i (\alpha )| \mbox{ for } i=0,1,2,3,6 .$$
Then $m_6 = 1$, $m_0+m_1+m_2+m_3 + 1 = 672 $ and 
$D2(\alpha )$ and $D4(\alpha )$ show that 
$$\ \ (\star ) \ \ m_0 = 1148 - 10 m_3, m_1 = 15 m_3 -1200, m_2 =723 - 6m_3 $$
so in particular $80 \leq m_3 \leq 114 $. 
Moreover for $z\in M_3(\alpha )$ also $\alpha -z \in M_3(\alpha )$ so 
$m_3$ is even. 
\\
\underline{Claim 1:} For all $\alpha \in Z $ we have $m_3 \neq 114$.
\\
Assume that there is some $\alpha \in Z$ with $m_3 =114$ and consider the
set
$$\overline{M}_3 := \{ \overline{z} := z-\frac{1}{2} \alpha \mid z\in M_3(\alpha ) \} 
\subset \alpha ^{\perp } .$$
Then $\overline{M}_3 = - \overline{M}_3$ is antipodal and 
for $\overline{z}_1,\overline{z}_2 \in \overline{M}_3$ we find
$$(\overline{z}_1,\overline{z}_2 ) = 
(z_1,z_2) - \frac{3}{2}  = \frac{\pm a}{2} \mbox{ with } 
a\in \{ 9,3,1 \} .$$
Choose some $y_0\in \sqrt{2} \overline{M}_3 = Y\disj -Y$  and let
$a_i (y_0) := |\{ y\in Y \mid (y,y_0) = \pm i \} | $ ($i=1,3$). 
Then 
$$ a_3(y_0) +a_1(y_0)  = 56 = |Y| -1 \mbox{ and } 
\sum _{y_0 \in Y} 9^2 + a_1(y_0) + 9 a_3(y_0) \geq 9^2|Y|\frac{57}{13} .$$
In particular there is some $y_0\in Y$ such that 
$a_3(y_0) \geq 27.2$, so $a_3(y_0) \geq 28$. 
Now $y_0 = \sqrt{2} (z_0 - \frac{1}{2} \alpha ) $ for some $z_0 \in M_3(\alpha )$ and $(y,y_0) = 3$ is equivalent to 
$(z,z_0) = 3$ hence there are at least 
28 vectors $z\in M_3(\alpha )$ with $(z,z_0)= 3$. 
This yields $28$ vectors $z-z_0 \in M_0(\alpha )$ contradicting the
fact that $m_0 = 1148 - 10m_3 = 8$ if $m_3 =114$. 
This proves Claim 1. Of course by interchanging the roles of the 
two lattices $\Lambda $ and $\Gamma $ we similarly obtain 
\\
\underline{Claim 1':} For all $\beta \in X $ 
$m_3'(\beta ) := |\{ x\in X \mid (x,\beta ) = \pm 3 \} | \leq 112 $.
\\
\underline{Claim 2:}  There is some $\alpha \in Z$ or some 
$\beta \in X$ with $m_3(\alpha ) = 114$ resp. $m_3'(\beta) =114 $. 
\\
To see this we use the set $T= \Min (\Gamma ) \cup \sqrt{6} \Min (Lambda ) 
= \pm (X\cup Z)$ defined above 
and the positivity property  (see Lemma \ref{positiv})
$$\sum _{t_1,t_2 \in T}  (t_1,t_2)^6 \geq 6^6 |T|^2 \frac{3\cdot 5}{14\cdot 16\cdot 18} $$
In particular there is some $t\in T$ such that 
$\sum _{t_1\in T}  (t_1,t)^6 \geq 6^6 |T| \frac{3\cdot 5}{14\cdot 16\cdot 18} $
and we may assume without loss of generality that $t = \alpha \in Z$.
We obtain
$n_i:= | \{ x\in X \mid (\alpha , x) = \pm i \sqrt{6} \} | $
as $n_2 = 3$ and $n_1= 2^2\cdot 3 \cdot 23 $. 
Using  $(\star )$ again we calculate
$$\sum _{t_1\in T} (t_1,t)^6 = 
6^6 + 2^6 6^3 n_2  + 6^3 n_1 + 3^6 m_3 + (723-6m_3) 2^6 (15m_3-1200)  
\geq  6^6 5 $$ 
which implies that $m_3 \geq 112.4 $ so $m_3 =114 $. 
\\
Of course Claim 2 contradicts Claim 1 or Claim 1'.
\eb
\section{Dual strongly perfect lattices of minimal type.} 

In this section we treat the 14-dimensional strongly perfect lattices of 
minimal type, where we focus on the dual strongly perfect lattices.
Nevertheless some of the results can be used for a later classification
of all strongly perfect lattices. 
The first two lemmata are of general interest and particularly useful 
in this section. 

\begin{lemma}\label{antipodal}
Let $Y = \Min (\Gamma )$ be a 4-design consisting of minimal vectors 
in a lattice $\Gamma \leq \R^n$ and $a:= \min (\Gamma ) = (y,y) $ for all $y\in Y$.
Assume that $(y_1,y_2) \neq 0$ for all $y_1,y_2 \in Y$.
%$(y_1,y_2) \not\in (a/2,a/4) $.
Fix $y_0\in Y$ and consider $N(y_0) := \{ y\in Y \mid (y,y_0) = a/2 \} $.
Then $|N(y_0)| \leq 2(n-1)$ and $|N(y_0)| $ is even.
\end{lemma}

\bew
For $y\in N(y_0)$ also $y_0-y \in N(y_0)$ so $|N(y_0)|$ is even.
To get the upper bound on $|N(y_0)|$ consider the projection
$$ \overline{N} := \{ \overline{y}:=y - \frac{1}{2} y_0 \mid y \in N(y_0 ) \}$$ 
of $N(y_0) $ to $y_0^{\perp }$.
For $y_1,y_2\in N(y_0)$  one has $(y_1,y_2) \neq a/2$ since 
otherwise $(y_0,y_1-y_2) = 0$. Therefore 
$$(\overline{y}_1,\overline{y}_2 ) =
(y_1,y_2) - a/4  \left\{ \begin{array}{ll} = 3a/4 & y_1= y_2 \\ \leq 0 & y_1\neq y_2  \end{array} \right. $$ 
hence distinct  elements  in $\overline{N}$ have non positive inner products.
Therefore $|\overline{N}| = |N(y_0)| \leq 2(n-1) $.
\eb

Let  $\Lambda $ be a 14-dimensional strongly perfect lattice of minimal type
scaled such that $\min (\Lambda ) = 1$.
Then $\min (\Lambda ^* ) = \frac{16}{3}$.
As above we choose $X\subset \Lambda _1$ such that 
$\Lambda _1 = X \disj -X $ and put $s:= |X| $.
Then by the bound for antipodal spherical codes given in
\cite{Anstreicher}, we have $s\leq 1746$. 
From equality $(D2)$ we know that $8 \frac{s}{21}  \in \Z $ and 
hence $$s=21 s_1 \mbox{ with } 5\leq  s_1 \leq 83  ~. $$

\begin{lemma}\label{bounds}
With the notation above let 
$$Z := \Min(\Lambda ) \cup \frac{\sqrt{3}}{4} \Min (\Lambda ^*) \subset 
{\cal S}^{(13)}(1) .$$ 
Then $Z$ is an antipodal kissing configuration and hence by 
\cite{Anstreicher} $|Z| \leq 2\cdot 1746 $. 
\end{lemma}

\bew
Let $x\in \Min (\Lambda )$ and $y\in \Min (\Lambda ^*)$.
Since $\Lambda $ is a strongly perfect lattice of minimal Type 
we have $|(x,y)| \leq 1$. 
In particular 
$$|(x,\frac{\sqrt{3}}{4} y) | \leq \frac{\sqrt{3}}{4} < \frac{1}{2} .$$
Since  $\Min(\Lambda )$ and $ \Min (\Lambda ^*) $ consist of 
minimal vectors of a lattice, we find $|(z_1,z_2)| \leq \frac{1}{2} $ 
for all $z_1 \neq \pm z_2 \in Z $.
\eb

\begin{bem}
Assume that there are
$\alpha , \beta \in \Lambda^* _{16/3}$ with
$0 \leq (\alpha ,\beta ) <  \frac{8}{3} $.
Then $\gamma := \alpha - \beta  \in \Gamma $ satisfies
$\frac{16}{3} <  (\gamma , \gamma ) \leq \frac{32}{3}  $.
Moreover for all $x\in X$ we have $|(x,\gamma ) | \leq 2$ hence 
$\gamma $ satisfies the conditions of Lemma \ref{linkomb}.
\end{bem}

Using the Theorem by Minkowski on the successive minima of a lattice we 
get:

\begin{bem}
If 
$|(\alpha ,\beta ) | = \frac{8}{3} $
for all $\alpha \neq \pm \beta \in \Lambda ^* _5$ 
then $\Lambda ^* _{16/3} = \sqrt{8/3} A_1 $ or $\Lambda ^* _{16/3} = \sqrt{8/3} A_2 $.
In particular 
there is some $\gamma \in \Lambda ^* $ with 
$\frac{16}{3} < (\gamma, \gamma ) < 8.2 <9 $.
\end{bem}

\begin{kor}
There is $\gamma \in \Lambda ^* $ with 
$\frac{16}{3} < r:= (\gamma , \gamma ) \leq \frac{32}{3} $ that satisfies the conditions of
Lemma \ref{linkomb}.
In particular the equation 
$$(\star ) \ \ \ \ |N_2(\gamma ) | =: n_2 = \frac{s_1 r (3r-16)}{2^7 }  $$
has a solution $(n_2,s_1,r)$ with natural numbers
$5\leq s_1 \leq 83, n_2 $ and some rational number 
$\frac{16}{3} < r \leq \frac{32}{3} $.
\end{kor}

\begin{rem}{\label{casemin}}
Searching for such solutions of $(\star ) $ by computer and using the
bound that $n_2 \leq \frac{r}{8-r}$ for $r<8$ (Lemma \ref{boundn2}) resp. 
$n_2 \leq 26$ for $r=8$ (Lemma \ref{boundr8}) 
we find the following 
possibilities: 
\begin{itemize}
\item[(a)]
$r=8$ or $r=32/3$ with $s_1 = 6,12,18,30,36,42$ 
\item[(b)]
$r=8$ with $s_1 = 10, 14,20,22 ,26, 28, 34, 38,44,46,52$. (\ref{only8})
\item[(c)]
$r=32/3$ with $s_1 = 9,15,21,33, 39,45, 51, 57,60,63,66,69,78$. (\ref{notc})
\item[(d)]
$r=8$ or $r=28/3$ with $s_1 = 8,16,40,56,80$. (\ref{notd})
\item[(e)]
$s_1=24$,  $r=20/3,8,28/3,32/3$ 
\item[(f)]
$s_1=25$,  $r=32/5,128/15, 48/5$. (\ref{s25})
\item[(g)]
$s_1=27$, $r=64/9, 80/9, 32/3 $. (\ref{s27})
\item[(h)]
$s_1=32$, $r=6,22/3,8,28/3,10 $ 
\item[(i)]
$s_1=48$, $r=8,28/3,32/3$. (\ref{s48})
\item[(j)]
$s_1=49$, $r=64/7, 208/21 $. (\ref{notjmnp})
\item[(k)]
$s_1=50$,  $r=8,128/15, 48/5,152/15$. (\ref{s50})
\item[(l)]
$s_1=54$, $r= 80/9, 88/9, 32/3 $. (\ref{s54})
\item[(m)]
$s_1=64$, $r=28/3,10 $. (\ref{notjmnp})
\item[(n)]
$s_1=72$, $r=28/3,32/3 $. (\ref{notjmnp})
\item[(o)]
$s_1=75$,  $r=128/15, 48/5,32/3$. (\ref{s75})
\item[(p)]
$s_1=81$,  $r=80/9,32/3$. (\ref{notjmnp})
\end{itemize}
In brackets are the references, where we exclude the corresponding 
case. The cases (a), (e), and (h) 
 are treated in the end of this section.
\end{rem}

\begin{lemma}\label{s25}
$s_1 \neq 25$.
\end{lemma}

\bew
Assume that $s_1 = 25$ and rescale such that $\min (\Gamma ) = 10$ and
$\min (\Lambda ) = 8/15 $ and let $X\disj -X = \Min (\Lambda )$. 
Then for all $\alpha ,\beta \in \Gamma $
$$ 
\begin{array}{lll}
\sum _{x\in X }(x,\alpha )^2 &= 20 (\alpha , \alpha ) & \in \Z \\
\sum _{x\in X }(x,\alpha )^4 &= 2 (\alpha , \alpha )^2 & \in \Z \\
\sum _{x\in X }(x,\alpha )^2 (x,\beta )^2 &= \frac{2}{3} (2(\alpha,\beta) ^2 + (\alpha , \alpha )(\beta ,\beta )) & \in \Z \\
\frac{1}{12} (D4-D2)(\alpha ) & = \frac{1}{6} (\alpha ,\alpha ) ((\alpha ,\alpha ) - 10 ) & \in \Z 
\end{array}
$$
which shows that $\Gamma $ is an even lattice with a sublattice
$$\Gamma ^{(t)} := \{ \gamma \in \Gamma \mid (\gamma ,\gamma ) \in 3 \Z \}$$
of index 3.
Then $\min (\Gamma ^{(t)} ) = 12$. 
Otherwise $\min (\Gamma ^{(t)}) \geq 18$ and hence 
$$9 \det (\Gamma ) = \det(\Gamma ^{(t)})  \geq (\frac{18}{\gamma _{14}})^{14} 
\geq 232242985896 =: u .$$
On the other hand 
$$\det(\Gamma ) \leq (\frac{15 \gamma _{14}}{8})^{14} \leq 10712486160 =: o.$$
Since $9o \leq u$ this is a contradiction.
In particular there is some $\beta \in \Gamma $ such that $(\beta,\beta ) =12$.
Then 
$$N_2(\beta ) = \{ x\in \Min(\Lambda ) \mid (x,\beta ) = 2 \} = \{ x_1,x_2,x_3,x_4\}$$ 
is of cardinality 4, 
$(x_i,x_j) = \frac{4}{15}$ for all $i\neq j$ 
 and $x_1+x_2+x_3+x_4 = \frac{2}{3} \beta $. 
This implies that $\xi := \frac{1}{3} \beta \in \Lambda $, 
$(\xi,\xi ) = \frac{4}{3}$ and 
$(x_1,\xi ) = \frac{2}{3} $. 
Therefore $x:= \xi - x_1 \in N_2(\beta )$, say $\xi - x_1 = x_2$. 
But then we have the relation $\frac{1}{3} \beta = x_1+x_2 $ implying
that  $(x_1,x_2) = \frac{2}{3} - \frac{8}{15} = \frac{2}{15}$,
a contradiction.
\eb

\begin{lemma}\label{notc}
If $r=32/3$ is the only possible norm of an element $\alpha \in \Lambda ^*$
with $16/3 < (\alpha , \alpha ) \leq 32/3$  
(so we are in case (c) of Remark \ref{casemin}) 
then $Y:=\Min (\Lambda ^*)$ is 
not a 4-design. 
\end{lemma}

\bew
In this case the inner products of $y_1,y_2 \in Y$ 
are $(y_1,y_2) \in \{ \pm 16/3, \pm 8/3 , 0 \} $ so 
$Y$ is a rescaled root system of rank 14. 
By \cite[Theorem 6.11]{Venkov} the only root systems that form 4-designs are 
$A_1, A_2, D_4, E_6, E_7$ and $E_8$. 
\eb

\begin{lemma}\label{only8}
If $r=8$ is the only possible norm of an element $\alpha \in \Lambda ^*$
with $16/3 < (\alpha , \alpha ) \leq 32/3$ 
(so we are in case (b) of Remark \ref{casemin})
 then $Y:=\Min (\Lambda ^*)$ is 
not a 4-design. 
\end{lemma}

\bew
Rescale $\Lambda ^*$ such that $(y,y) = 4$ for all $y\in Y$. 
Then
the inner products of $y_1,y_2 \in Y$ 
are $(y_1,y_2) \in \{ \pm 4, \pm 2 , \pm 1 \} $.
Assume that $Y$ is a 4-design and let $a:=|Y|/2$.
For $y_0\in Y$ and $i=1,2$ let 
$a_i:= | \{  y\in Y \mid (y_0,y) = i \} |$.
Then the 4-design conditions yield that 
$$ 1+a_1+a_2 = a, 4^2+a_1+2^2a_2 = 4^2 a/14, 
 4^4+a_1+2^4a_2 = 4^4 3 a/(14\cdot 16) $$ 
having the unique solution 
$$ a_2 = 0, a = 105, a_1 = 104 .$$
So in this case $Y$ is an antipodal
 tight 5-design, but  no such 
tight design exists in dimension 14 (see for instance \cite{tight}). 
\eb

\begin{lemma}\label{notd}
If $r=8$ and $r=28/3$ are the only possible norms of elements
 $\alpha \in \Lambda ^*$
with $16/3 < (\alpha , \alpha ) \leq 32/3$  
(so we are in case (d) of Remark \ref{casemin})
then $Y:=\Min (\Lambda ^*)$ is 
not a 4-design. 
\end{lemma}

\bew
Rescale $\Lambda ^*$ such that $(y,y) = 8$ for all $y\in Y$. 
Then
the inner products of $y_1,y_2 \in Y$ 
are $(y_1,y_2) \in \{ \pm 8, \pm 4 , \pm 2 , \pm 1 \} $.
Assume that $Y$ is a 4-design and let $a:=|Y|/2$.
For $y_0\in Y$ and $i=1,2,4$ let 
$a_i:= | \{  y\in Y \mid (y_0,y) = i \} |$.
Then by the 4-design conditions 
$$a = 105 + 5 a_4, a_1= \frac{64}{21} a_4, a_2 = 104 + \frac{20}{21} a_4 ; $$
in particular $a_4$ is divisible by 21.
Moreover  $a_4 > 0$ by Lemma \ref{only8} and 
$a_4$ is even and $a_4 \leq 26$ by Lemma \ref{antipodal} which is 
a contradiction.
\eb

\begin{lemma}\label{notjmnp}
In the cases (j), (m), (n), (p) of Remark \ref{casemin} 
$\Min (\Gamma )$ is not a 4-design.
\end{lemma}

\bew
In these cases only three inner products 
$i,j,8/3$ 
(up to sign) between distinct $y_1,y_2\in Y$ are possible. 
For a fixed $y_0 \in Y$ 
let 
$$\begin{array}{ll} 
a:= |Y|/2, &  b:= \frac{1}{2} |\{ y\in Y  \mid (y,y_0) = \pm i \} | \\
c:= \frac{1}{2} |\{ y\in Y  \mid (y,y_0) = \pm j \} |,  &
d:= \frac{1}{2} |\{ y\in Y  \mid (y,y_0) = \pm 8/3 \} |.
\end{array} $$
Assuming that $Y$ is a 4-design
one obtains the system 
$$ 1+b+c+d = a, \ 
 d \frac{64}{9} + c j^2+ b i^2 =  (\frac{a}{14} -1) \frac{16^2}{3^2} ,\ 
 d \frac{8^4}{3^4} + c j^4+ b i^4 =  (\frac{3a}{14\cdot 16} -1) \frac{16^4}{3^4} $$
of which the solution depends on one parameter. 
In the different situations we find: 
\begin{itemize}
\item[(j)]
$s_1=49$, $i =16/21,j=  8/21 $
$$a=(24960 + 1440d)/299 , b=(6825+112d)/23, c= -(64064 + 315d)/299$$
so $c$ is negative here, a contradiction.
\item[(m)]
$s_1=64$, $i=2/3,j=1/3 $
$$a=(37585+2205 d)/461, b=(176800+3024 d)/461, c = - (139776+1280 d)/461 $$ 
so $c$ is negative here, a contradiction.
\item[(n)]
$s_1=72$, $i=2/3,j=0 $
$$a=(1764+105d)/22, b=(3328+64d)/11, c=-(4914 +45 d)/22 $$
so $c$ is negative here, a contradiction.
\item[(p)]
$s_1=81$,  $i=8/9,j=0$.
$$a=1/23(1960+112d), b=1/23(4212 + 81 d), c=1/23 (-2275+8d) $$
For $c$ to be a nonnegative integer, $d$ has to be of the form 
$293 + 23 x $. Since $d$ is even, $x=1,3,5 \ldots $ is odd. 
In particular $a  \geq 1624 $. 
Since $s= \frac{1}{2} |\Min (\Lambda ) | = 21 \cdot 81 =1701 $ 
this is a contradiction to Lemma \ref{bounds}.
\end{itemize}
\eb

From now on we assume that $\Lambda $ and $\Lambda ^*$ are 
both strongly perfect. 
In particular $Y \disj -Y :=\Min (\Lambda ^*)$ is also a $4$-design. 
We put $|Y| = 21 t_1$ and $|X| = 21 s_1 $ 
(where $X \disj -X = \Min (\Lambda )$). 
By the discussion above we know that both, $s_1$ and $t_1$, belong to the set
$$S:= \{ 6,12,18,24,27,30,32,36,42,48,50,54,75 \} $$
and $s_1+t_1 \leq 83 $ by Lemma \ref{bounds}.
Moreover Remark \ref{casemin} gives all possible inner products 
$$i:=|(y_1,y_2)| = (32/3 - r)/2 $$
 for $y_1,y_2 \in Y$ according to the different values of $s_1$.
Since also $Y$ is a 4-design, we obtain a system of equations 
for the values $a_i:=|\{ y \in Y \mid |(y,y_0)|  =  i \} | $ 
for fixed $y_0\in Y$.
%$$ \star \  \begin{array}{ll} 
%1+ \sum _{i} a_i  & = 21 t_1  \\ 
%\frac{16}{3}^2 + \sum _{i} i^2 a_i  & = \frac{2^7}{3} t_1 \\
%\frac{16}{3}^4 + \sum _{i} i^4 a_i  & = \frac{2^{11}}{3^2} t_1 
%\end{array}
%$$
The values $a_i$ are nonnegative integers. 
Going through the different possibilities we find the next lemmata:

\begin{lemma} \label{s75}
$s_1 \neq 75$. 
\end{lemma}

\bew
Assume that $s_1 = 75$. 
Rescale the lattice such that $\min (\Lambda )= \frac{8}{15}$.
Then $\min (\Gamma ) = 10 $ and for $y_1,y_2\in \Gamma _{10}$ we have
$$|(y_1,y_2)| \in \{10,5,2,1,0 \} .$$
Fix $y_0 \in Y$ and let $a_i := | \{ y\in Y \mid (y,y_0) = \pm i  \} |$.
Then the solution of the system of equations 
$$\begin{array}{ll} 
1+a_5+a_2+a_1+a_0  & = 21 t_1 \\
10^2+5^2a_5+2^2a_2+a_1  & = 150 t_1 \\
10^4+5^4a_5+2^4a_2+a_1  & = 5^4 3^2 t_1/2 
\end{array}
$$
has the following solution depending on the two parameters $a_5$ 
and $t_1 = 8 t_2$:
$$a_2 = 1775 t_2 - 825 -50 a_5 , \ 
a_1 = 3200 +175 a_5 - 5900 t_2 , \
a_0 = 4293 t_2 - 2376 - 126 a_5. $$
In particular $t_1$ is divisible by 8,
which leaves the possibilities $t_1=24,32,48$ contradicting the fact that
$s_1+t_1 \leq 83$ hence $t_1 \leq 8$.
\eb

\begin{lemma}\label{s54}
$s_1 \neq 54 $.
\end{lemma}

\bew
Rescale such that $\min (\Gamma ) = 12, \min (\Lambda ) = \frac{4}{9} $.
Then the inner products between minimal vectors of $\Gamma $ are 
$$\{ \pm 12, \pm 6 , \pm 2,\pm 1, 0 \} .$$
The solution of the respective system of equations for the 
$a_i$ depends on two parameters $a_6$ and $t_1$:
 $$a_2 = -1716 - 105 a_6 + 468 t_1, a_0 = -5005 - 280 a_6 + 1209 t_1,
    a_1 = 6720 + 384 a_6 - 1656 t_1 .$$
In particular $a_0 \geq 0 $ and $a_1\geq 0$ yields 
$$ (1656t_1-6720)/384 \leq a_6 \leq (1209t_1 -5005)/280  $$ 
which implies that $t_1 \geq 70$, a contradiction. 
\eb

\begin{lemma}\label{s27}
$s_1\neq 27$.
\end{lemma}

\bew
Rescale such that $\min (\Gamma ) = 6, \min (\Lambda ) = \frac{8}{9} $.
Then $\Gamma $ is an even lattice
and the inner products between minimal vectors of $\Gamma $ are
$$\{ \pm 6, \pm 3,  \pm 2, \pm 1 , 0 \} .$$
There is some pair $\alpha _1,\alpha _2 \in \Gamma _6$ such that
$(\alpha _1,\alpha _2) = 2$, otherwise the fact that 
$\Gamma _6 = \pm Y$ is a  4-design, in particular 
the equations $D2(\alpha )$ and $D4(\alpha )$ would imply that for any fixed 
$\alpha\in \Gamma _6$ 
$$|Y| = 1470 + 14 | Y \cap \alpha ^{\perp } | 
\geq 1470 $$ contradicting Lemma \ref{bounds}.
So choose $\alpha _1,\alpha _2 \in \Gamma _6$ such that $(\alpha_1,\alpha_2) = -2$ and put $\beta := \alpha _1 + \alpha _2 \in \Gamma _8$.
From the fact that $ \Lambda _{8/9} $ is a 4-design, we get that 
$$N_2(\beta ) : = \{ x\in \Lambda _{8/9} \mid (x,\beta ) =2 \} 
= \{ x_1,\ldots , x _8 \} $$
such that $\sum _{i=1}^8 x_i = 2 \beta $ and 
$(x_i,x_j) = \frac{4}{9} $ for all $i\neq j$.
Hence $N_2(\beta )$ generates a rescaled version of the root lattice
$A_8$ in $\Lambda $ and the residue classes 
$x_i + \Gamma \in \Lambda / \Gamma $ generate a subgroup of 
$\Lambda /\Gamma $ of order at least $9^7$.
Since $\min (\Lambda ) = \frac{8}{9}$ we get 
$$\det (\Gamma ) \leq (9\cdot 2.776/8)^{14} < 2\cdot  9^7 $$
and therefore $\det (\Gamma ) = 9^7 $. 
But then the lattice 
$$\Delta := \langle 3 x_1,\ldots ,3 x_7, \Gamma \rangle \leq \Lambda $$ 
is an even unimodular sublattice of $\Lambda $ of rank 14, which is
a contradiction.
\eb

We now rescale the lattice such that $\min (\Gamma ) = 4$, $\min (\Lambda ) = 4/3$.
Then for all $\alpha ,\beta  \in \Gamma $
$$
\begin{array}{lll}
\sum _{x\in X} (x,\alpha )^2  & = 2 s_1 (\alpha , \alpha )  & \in \Z  \\
\sum _{x\in X} (x,\alpha )^4  & =  s_1/2 (\alpha , \alpha )  & \in \Z  \\
\sum _{x\in X} (x,\alpha )^2 (x,\beta )^2  & =  s_1/6 (2 (\alpha ,\beta )^2 + (\alpha , \alpha ) (\beta ,\beta ) ) & \in \Z  \\
\frac{1}{12} 
\sum _{x\in X} (x,\alpha )^4 -(x,\alpha)^2 & =  s_1/8 (\alpha , \alpha ) ((\alpha , \alpha ) - 4)  & \in \Z  
\end{array}
$$

In particular, if $s_1$ is not a multiple of 3 then 
$\Gamma $ is 3-integral and the set 
$$\Gamma ^{(t)} := \{ \alpha \in \Gamma \mid (\alpha ,\alpha ) \in 3 \Z _{(3)} 
\} $$
is a sublattice of $\Gamma $ of index 3.

\begin{lemma}\label{m3}
$s_1 t_1$ is a multiple of 3.
\end{lemma}

\bew
Assume that both, $s_1$ and $t_1$ are not divisible by $3$.
Then rescaled to $\min (\Gamma ) = 4, \min(\Lambda )= 4/3$ the lattice
$\Gamma $ is $3$-integral and also $\sqrt{3}\Lambda $ is 3-integral and
both lattices have a sublattice of index 3 in which all inner products are 
divisible by 3. 
In particular the 3-part of the determinant of $\Gamma $ is 
$3^a$ with $a \geq 12$ and the one of $\sqrt{3}\Lambda $ is 
$3^{14-a}$ again with $14-a \geq 12$, which is a contradiction. 
\eb

%\begin{lemma}\label{s25}
%If $s_1=25$ then $t_1\in \{ 24,48 \}$.
%\end{lemma}
%
%\bew
%If $s_1=25$ then rescaled to $\min (\Gamma )=10$, the inner products
%between minimal vectors in $\Gamma $ are 
%$ \pm 10, \pm 5 , \pm 4, \pm 2, \pm 1 $.
%Solving the corresponding system of equations we find that 
%$21 t_1 = (  14784 + 280 a4  + 784 a5 )/159 $ 
%and so $t_1$ is a multiple of $8$. 
%By Lemma \ref{m3} it is also divisible by 3 and hence the Lemma follows.
%\eb

\begin{lemma}\label{s50}
$s_1 \neq 50 $. 
\end{lemma}

\bew
Rescale such that $\min (\Gamma ) = 20, \min (\Lambda ) = \frac{4}{15} $.
Then the inner products between minimal vectors of $\Gamma $ are
$$\{ \pm 20, \pm 10,  \pm 5, \pm 4 , \pm 2,\pm 1 \} .$$
First assume that there is some $\beta \in \Gamma _{30}$ 
(this is certainly the case if there are $\alpha _1,\alpha _2 \in \Gamma _{20}$
with $(\alpha _1, \alpha _2 ) = 5 $). 
Then 
$$N_2(\beta ) := \{ x\in \Lambda _{4/15} \mid (x,\beta ) = 2 \} = 
\{ x_1,\ldots , x_{25} \} $$ 
and for $x\in N_2(\beta )$, $\overline{x} := x-\frac{1}{15} \beta $ is 
perpendicular to $\beta $.
Moreover for $x_1,x_2 \in N_2(\beta )$ 
$$(\overline{x}_1,\overline{x}_2) = (x_1,x_2) - \frac{2}{15}  
\left\{ \begin{array}{ll}  = 2/15 & \mbox{ if } x_1 = x_2  \\
\leq 0 & \mbox{ if } x_1 \neq x_2  \end{array} \right. .$$
So the set $\overline{N}_2$ consists of 25 vectors of equal length in 
a 13 dimensional space having non positive inner products.
Therefore there are two elements $x_1,x_2 \in N_2(\beta )$ with
$\overline{x}_1+\overline{x}_2 = 0$ i.e. $$x_1+x_2 = \frac{2}{15} \beta .$$
Hence $\beta \in 15 \Gamma ^* $ and therefore 
$$(\alpha , \beta ) \in \{ 0 ,\pm 15 \}  \mbox{ for all } \alpha \in \Gamma _{20 }. $$
Since $\Gamma _{20}$ is a 4-design of cardinality $42 t_1$ the numbers
$n_i:= | \{ \alpha \in \Gamma _{20} \mid (\alpha ,\beta ) = \pm i \} |$
satisfy
$$n_0+n_{15} = 42 t_1, 15^2 n_{15} = 1800 t_1, 15^4 n_{15} =  202500 t_1 $$
which has $0$ as only solution. 
So this case is impossible.
\\
Now we fix $\alpha \in \Gamma _{20} = Y \disj - Y $ and put 
$$a_i := | \{ y \in Y \mid (y,\alpha ) = \pm i \} | \ \mbox{ for }
 i= 20, 10,  5,  4 ,  2, 1 .$$
Then we have $a_{20} =1$ and $a_5 = 0 $ by the above consideration.
Moreover $|Y| = 21 t_1 $ with $t_1\leq 32$. 
The solutions of the 4-design equations for $Y$ in the $a_i$ depend on
two parameters $a_{10}$ and $t_1$ and one finds
$$5a_4 =  -4389 + 1169 t_1 - 264 a_{10} ,\
 5a_1 = -16896 + 3816 t_1 - 896 a_{10} ,\ 
 a_2= 4256 - 976 t_1 + 231 a_{10} .$$
Going through the different possibilities for $t_1\leq 32$ the only 
even nonnegative integer  $a_{10}$ such that $a_1,a_2,a_4$ are all
nonnegative integers is 
$$ t_1 =27, a_{10}= 96, a_1= 24 , a_2=80, a_4= 366 .$$
But this case is impossible by Lemma \ref{s27}.
\eb

\begin{lemma}\label{s48}
$s_1 \neq 48$.
\end{lemma}

\bew
Rescale $\Gamma $ such that $\min (\Gamma ) =8$ and $\min (\Lambda ) = 2/3$.
Then $\Gamma $ is an even lattice and the inner products 
between minimal vectors in $\Gamma $ are
$$\{ \pm 8, \pm 4, \pm 2, \pm 1, 0\}.$$
Assume that there is some $\alpha \in \Gamma _8$ such that 
$(\alpha,\alpha') \neq \pm 2$ for all $\alpha ' \in \Gamma _8 $.
Since $\Gamma _8$ is a 4-design, the equations $D4(\alpha )$ and 
$D2(\alpha )$ imply that 
$$|\Gamma _8 \cap \alpha ^{\perp }| = -2 (189 + 9 t_1) < 0 $$
which is a contradiction.

Choose such a pair $(\alpha _1 , \alpha _2) \in \Gamma_8^2$ such that 
$(\alpha _1 , \alpha _2) = -2$ and put 
$\beta := \alpha _1 + \alpha _2  \in \Gamma _{12} $.
Then $N_2(\beta ) =\{ x\in \Lambda _{2/3} \mid (x,\beta ) = 2 \}$ 
has cardinality 24.
Put $$\overline{N}_2 :=\{ x - \frac{1}{6} \beta \mid x\in N_2(\beta ) \}.$$
Then $\overline{N}_2 \subset \langle \alpha _1 ,\alpha _2 \rangle^{\perp }$ 
and any two distinct vectors  in 
$\overline{N}_2$ have non positive inner product. 
In particular we may apply \cite[Lemma 2.10]{dim12} to write 
$\overline{N}_2 = \disj _{i=1}^k E_i $ as a union of disjoint indecomposable 
components 
$E_i$ (i.e. $\sum _{x\in E_i} x  = 0 , E_i \perp E_j $ if $i\neq j $)
 with $$k= |\overline{N} _2 | - \dim (\langle \overline{N}_2 \rangle ) 
\geq 24-12 = 12 .$$
So $k=12 $ and $E_i = \{ \overline{x}_i ,-\overline{x}_i \} $ for all $i$.
In particular there are $x_1,x_{13} \in N_2(\beta )$ 
such that $x_1+x_{13} = \frac{1}{3} \beta $.
So $\beta \in 3 \Gamma ^*$ and hence
$(\beta , \alpha ) \in \{ 0,\pm 3, \pm 6 \} $
for all $\alpha \in \Gamma _8$.
For $i=0,3,6$ put $n_i:= |\{ \alpha \in Y \mid (\alpha , \beta ) = \pm i \} | $. 
Then from the fact that $\Gamma _8 = Y \disj -Y $ is a 4-design of cardinality $42 t_1$ 
one obtains a unique solution 
$$n_0 = 9 t_1, \ n_3 = 32/3 t_1 , \ n_6 = 4/3 t_1 .$$
Now consider the set $$M_6(\beta ) := \{ \alpha \in \Gamma _8 \mid 
(\alpha , \beta ) = 6 \}  = \{ \alpha _1,\ldots , \alpha _{n_6} \} $$
 of cardinality $n_6 = 4/3 t_1 $ and choose
$\alpha \in M_6(\beta )$.
The $x_i \in N_2(\beta ) = \{ x\in \Lambda _{2/3} \mid (x,\beta ) = 2 \} $ 
satisfy $\sum _{i=1}^{n_2} x_i = \frac{1}{6} n_2 \beta $ and therefore
$\sum _{i=1}^{n_2} (\alpha ,x_i ) = n_2 $ hence $(\alpha , x_i) = 1$ for 
all $x_i \in N_2(\beta )$.
Hence the orthogonal projection $\overline{\alpha } := \alpha - \frac{1}{2} \beta \in \beta ^{\perp }$ 
and $\overline{x} := x - \frac{1}{6} \beta \in \overline{N} _2 $  are perpendicular for all $x\in N_2(\beta )$ since
$$
(\alpha - \frac{1}{2} \beta , x - \frac{1}{6} \beta ) = 
(\alpha - \frac{1}{2} \beta , x ) = 
(\alpha , x) - \frac{1}{2}(\beta ,x) = 1 - 1 = 0 .$$
Therefore $\overline{N}_2 \perp \overline{M}_6 \perp \langle \beta \rangle $ 
and 
$\dim (\overline{N}_2) = 12$ implies that  $\dim (\overline{M}_6(\beta )) = 1$. 
Since the elements in $\overline{M}_6(\beta )$ have all equal norm, this 
leaves $4/3t_1 = n_6= |M_6(\beta ) | = 2 $ which is a contradiction.
\eb

\begin{lemma}\label{s32}
If $s_1 = 32 $ then $t_1=6$.
\end{lemma}

\bew
Assume that $s_1=32$ and scale such that $\min (\Lambda ) = 1/3$, 
$\min (\Gamma ) = 16$. 
Then $\Gamma $ is an even lattice containing the sublattice 
$$\Gamma ^{(t)} := \{ \gamma \in \Gamma \mid (\gamma,\gamma ) \in 3\Z \} $$
of index 3. 
By Lemma \ref{m3} and Lemma \ref{s27} we know that 
$t_1=6t_2$ is a multiple of 6. 
Moreover the inner products 
$(\alpha _1 ,\alpha _2)$ for two elements $\alpha _i \in \Gamma _{16}$ 
are $\pm 1,\pm 2, \pm 4, \pm 5 , \pm 7, \pm 8 , \pm 16 $. 
Fix $\alpha _0 \in \Gamma _{16} = \pm Y $  and let 
$a_j := |\{ \alpha \in Y \mid (\alpha , \alpha _0 ) = j \} | $.
Then $a_8$ is even.
\\
Since $\Gamma _{16} $ is a 4-design of cardinality $21\cdot 6 \cdot t_2 $ 
we may solve the corresponding system of equations to obtain 
$$a_4 = 539 + 65/3 a_8 + 2/3 a_1 + 10 a_7 -522 t_2 .$$
In particular $a_4$ is odd and hence nonzero. 
\\
So there are $\alpha _0 , \alpha _0 ' \in \Gamma _{16} $ such that 
$\beta := \alpha _0 +\alpha _0 '  \in \Gamma _{24} $ has norm 24. 
Consider the set 
$$N_2(\beta ) := \{ x \in \Min (\Lambda ) \mid (x,\beta ) = 2 \} .$$
Note that for all $x\in N_2(\beta )$ the inner products $(x,\alpha _0 ) = (x,\alpha  _0' ) = 1 $.
Then $|N_2(\beta )| = 16$ and $\sum _{x\in N_2(\beta ) }x = 4/3 \beta $. 
The set $$\overline{N}_2 := \{ \overline{x} = x-\frac{1}{12} \beta \mid 
x\in N_2(\beta ) \} \subset \langle \alpha _0,\alpha _0 ' \rangle  ^{\perp }$$ is a set of vectors of equal norm
living in a 12-dimensional space such that distinct vectors have non positive
inner product and satisfying the relation $\sum _{\overline{x} \in \overline{N}_2} \overline{x} = 0 $. 
Write $$\overline{N}_2 = E_1 \disj \ldots \disj E_k $$  
as a disjoint union of indecomposable sets. 
Then by \cite[Lemma 2.10]{dim12} 
$\dim (\overline{N}_2) = 16 - k $ and hence $k\geq 4$. 
Moreover for all $i$ $$\sum _{x\in E_i} x = \frac{|E_i|}{12} \beta .$$
\\
Distinguish two cases:
\\
(a) There is some $i$ such that $|E_i|$ is not a multiple of 4. 
\\
(b) $k=4$ and  $|E_i| = 4 $ for all $i = 1,\ldots , 4 $.
\\
In case (a) the vector $\beta \in 6\Lambda \cap \Gamma $.
In particular $(\alpha , \beta ) \in \{ 0,\pm 6, \pm 12 \} $ for all 
$\alpha \in \Gamma _{16} = Y \disj -Y $. 
Since $\Gamma _{16}$ is a 4-design of cardinality $6\cdot 21 \cdot t_2$, the 
cardinalities 
$$m_j := |\{ \alpha \in Y \mid (\beta , \alpha ) = \pm j \} |$$
can be calculated as 
$$ m_0 = 54 t_2, m_6 = 64 t_2 , m_{12} = 8 t_2 . $$
Now consider the set 
$$M_{12}(\beta ):= \{ \alpha \in \Gamma _{16} \mid (\alpha ,\beta) = 12 \} $$
and its orthogonal projection 
$$\overline{M}_{12} := \{ \overline{\alpha } := \alpha - \frac{1}{2} \beta \mid \alpha \in M_{12}(\beta ) \} $$
into $\beta ^{\perp }$. 
Note that for $\alpha \in M_{12}(\beta )$ we have
$$16 = (\alpha , \frac{4}{3} \beta ) = \sum _{x\in N_2(\beta )} (\alpha , x) $$
and hence $(\alpha , x) = 1$ for all $x\in N_2(\beta )$. 
In particular 
$$ \overline{M}_{12} \perp \overline{N}_2 \perp \langle \beta \rangle .$$
Moreover $\overline{M}_{12} = - \overline{M}_{12} $ is an antipodal set.
For $\alpha _1,\alpha _2 \in M_{12}(\beta) $ we calculate 
$$(\overline{\alpha }_1, \overline{\alpha }_2 ) = 
(\alpha _1,\alpha _2) - 6 \in \{ 10,2,1,-1,-2,-10 \} .$$
Fix some $z _0 \in \overline{M}_{12}$ and let 
$n_i := | \{ z \in  \overline{M}_{12} \mid (z,z_0) = i \} |$ $(i=10,2,1)$.
Then $\dim ( \overline{M}_{12} ) \leq 5$  and $|\overline{M}_{12} | = 8 t_2 $
and hence 
$$10^2 + n_1 + 2^2 n_2 \geq \frac{4 t_2}{5} 10^2 = 80 t_2 ,\ n_1+n_2 = 4t_2 -1 $$
which yields $3 n_2 \geq 76 t_2 -99 $.
Now $n_2 \leq 4t_2-1$ so $4t_2-1 \geq 76 t_2 - 99$ yielding the inequality 
$t_2 \leq \frac{98}{72} < 2 $. 
So the only solution is $t_2=1$. 
\\
Assume now that  we are in case (b). 
Let $E_1 = \{ x_1,x_2,x_3,x_4 \} $ with $\sum _{i=1}^4 x_i = \frac{1}{3} \beta $.
Assume that $t_1\neq 24$. then 
$t_1 \in A$ and $$(x_i,x_j) \in \{ 0 ,\pm \frac{1}{12},\pm \frac{1}{6} , \pm \frac{1}{3} \} .$$
If $(x_i,x_j) = 0$, then $\frac{1}{6} \beta = x_i + x_j $ and $E_1$ is not
indecomposable. Hence $(x_i,x_j) \neq 0$ for all $x_i,x_j \in E_1$.
Moreover 
$$\frac{1}{3} = (\frac{1}{3} \beta ,x_1 ) - (x_1,x_1) = \sum _{i=2}^4 (x_i,x_1) 
= \frac{1}{6} + \frac{1}{12} + \frac{1}{12} $$ 
has only this solution. 
So the Gram matrix is 
$$A:=((x_i,x_j)_{i,j=1}^4) = \frac{1}{12} \left( \begin{array}{cccc} 
4 & 1 & 1 & 2 \\ 
1 & 4 & 2 & 1 \\ 1 & 2 & 4 & 1 \\ 2 & 1 & 1 & 4 \end{array} \right).$$
The vectors in $N_2(\beta )$ together with $\alpha _0$ 
(from above with $\beta  = \alpha _0 + \alpha _0' $) 
hence generate a lattice in $\Q \Gamma $ with Gram matrix 
$$ M = \frac{1}{12}  \left( \begin{array}{cccccccccccccc} 
4&1&1&2&2&2&2&2&2&2&2&2&2&12\\
1&4&2&1&2&2&2&2&2&2&2&2&2&12\\
1&2&4&1&2&2&2&2&2&2&2&2&2&12\\
2&1&1&4&2&2&2&2&2&2&2&2&2&12\\
2&2&2&2&4&1&1&2&2&2&2&2&2&12\\
2&2&2&2&1&4&2&2&2&2&2&2&2&12\\
2&2&2&2&1&2&4&2&2&2&2&2&2&12\\
2&2&2&2&2&2&2&4&1&1&2&2&2&12\\
2&2&2&2&2&2&2&1&4&2&2&2&2&12\\
2&2&2&2&2&2&2&1&2&4&2&2&2&12\\
2&2&2&2&2&2&2&2&2&2&4&1&1&12\\
2&2&2&2&2&2&2&2&2&2&1&4&2&12\\
2&2&2&2&2&2&2&2&2&2&1&2&4&12\\
12&12&12&12&12&12&12&12&12&12&12&12&12&192 \end{array} \right) .$$
The determinant of $M$ is in $15 (\Q^*)^2 $. 
In particular this means that the determinant of $\Gamma $ is a multiple
of 5 which is impossible since the exponent of $\Lambda/\Gamma $ 
divides 12. 

For $t_1=24$ we can do the same considerations. %to see that 
%$\dim (\langle \overline{N}_2 \rangle ) = 16-4 = 12 $ and hence 
%$\dim (\langle \overline{M}_{12} \rangle )  = 1$, so the representation 
%of $\beta = \alpha _0 + \alpha _0'$ as a sum of two elements of norm 16 
%is unique (there are no other $\alpha \in \Gamma _{16} $ that have 
%$(\alpha ,\beta ) = 12 $).
But here we get more possibilities for the inner products of the elements
in one component $E_i$. 
The possible Gram matrices are   
$$A, B:=\frac{1}{24} \left( \begin{array}{cccc} 
8 & 1 & 3 & 4 \\ 
1 & 8 & 4 & 3 \\ 3 & 4 & 8 & 1 \\ 4 & 3 & 1 & 8 \end{array} \right), \mbox{ or }
C:=\frac{1}{24} \left( \begin{array}{cccc} 
8 & 3 & 3 & 2 \\ 
3 & 8 & 2 & 3 \\ 3 & 2 & 8 & 3 \\ 2 & 3 & 3 & 8 \end{array} \right).
$$
So there are 15 possibilities for the Gram matrix $M$.
All these matrices have determinant $5 d$ for some rational $d$ with
denominator and numerator not divisible by 5. 
\eb

The following table lists the possible values of 
$s_1 = \frac{1}{42} |\Min (\Lambda )| $ 
together with a value 
$m:=\Min (\Lambda )$ for which $\Gamma = \Lambda ^*$ is an even lattice
of minimum $r=\frac{16}{3m} $. 
For each possible $s_1$ we list the possibilities for  
$t_1 = \frac{1}{42} |\Min (\Gamma )| $ and give numbers $e= e(s_1,t_1)$
 such that $e \Lambda \subset \Gamma $. 
$A$ stands for the set $\{6,12,18,30,36,42 \}$.
$$
\begin{array}{|c|c|c|c|c|c|}
\hline
&  & t_1 &  A & 24 &  32   \\
r & m & s_1 &  & &    \\
\hline
4 & 4/3 & A & 3 & 6 &  12 (s_1=6)   \\ 
\hline
8 & 2/3 & 24 & 6 & 12 &   (\ref{s32})  \\ 
\hline
16 & 1/3 & 32 &  12 (t_1=6) & (\ref{s32})  & (\ref{m3})    \\
\hline
\end{array}
$$

There are two different strategies with which we treat the remaining
cases: 
If $\exp (\Lambda / \Gamma )$ is 3 or 6 then we explicitly construct all
possible lattices as sublattices of one of the maximal even lattices 
$M$ given in Section \ref{maxeven}. 
If this exponent is 12, then this explicit calculation becomes too 
involved and there are too many lattices to be constructed. 
In this case we use the theory of modular forms as explained 
in the proof of Lemma \ref{beispmodular} to exclude the possible
theta series. 

The most interesting but also easy 
 case is of course that both $s_1 $ and $t_1$ lie in $A$, 
since the unique dual strongly perfect lattice $[\pm G_2(3)]_{14}$ 
of dimension 14 satisfies
$s_1=t_1=18$. 

\begin{satz}\label{g23}
Assume that $\Gamma $ is an even dual strongly perfect lattice
of dimension 14 of minimal type such that its 
discriminant group $\Lambda / \Gamma $ and also $\Gamma / 3\Lambda $
are 3-elementary. 
Then $\Gamma = [\pm G_2(3)]_{14}$.
\end{satz}

\bew
The lattice $\Gamma $ is a sublattice of one of the two maximal even
lattices $M$ of determinant 3 given in Section \ref{maxeven}. 
It has minimum 4 and its dual lattice has minimum $\frac{4}{3}$.
Successively constructing all sublattices $L$ of $M$ of 3-power index 
such that $L^*/L$ is of exponent 3 and $\min (L^*) \geq 4/3$ 
we find up to isometry a unique such $L$ that is of minimum $4$.
This is the extremal 3-modular lattice with kissing number 
$756 = 2\cdot 21\cdot 18$ and automorphism group  $G_2(3)\times C_2$
as stated in the theorem.
\eb

\begin{lemma}
If $s_1=24$ then $t_1 \not\in A$. 
\end{lemma}

\bew
Assume that $s(\Lambda ) = 21 \cdot t_1$ for some $t_1$ in $A$.
Then $\Gamma = \Lambda ^*$ rescaled to $\min (\Gamma ) = 4$ is an 
even lattices and $\min (\Lambda ) = \frac{4}{3}$.
If $s(\Gamma ) =21 \cdot 24 $, then $\sqrt{6}\Lambda $ is even,
hence the exponent of $\Lambda / \Gamma $ is 6, and 
$\Gamma $ is contained in one of the maximal even lattices 
$M$ from Section \ref{maxeven} for which $\sqrt{6} M^*$ is again even.
These are the two lattices in the genus of $E_6\perp E_8$ and 
the 6 lattices in the genus of $A_2 \perp D_{12}$.
From the latter genus only 2 lattices $M$ satisfy $\min (M^*) \geq \frac{4}{3}$.
Going through all possible sublattices $L$ of the relevant 4 maximal 
even lattices for which $\sqrt{6} L^*$ is even and has minimum 
$\geq 8 $ ($=6 \cdot \frac{4}{3} $) one finds only the lattice from
Theorem \ref{g23}.
\eb

\begin{lemma}
If $s_1=24$ then $t_1 \neq 24$. 
\end{lemma}

\bew
Assume that $t_1=s_1=24$ and rescale such that 
$\min(\Gamma ) = 8$, $\Gamma = \Lambda ^*$ with $\min (\Lambda ) = \frac{2}{3}$.
Then $\Gamma $ is even and so is $\sqrt{12} \Lambda $. 
We proceed as in the proof of Lemma \ref{beispmodular}. 
From Lemma \ref{boundsfromgamma} we find that 
$$ 2725130 \leq \det(\Gamma ) = 2^a 3^b \leq 471140123 .$$
Moreover $1\leq a\leq 27$ and $1\leq b \leq 13$ and the 
 Jordan decomposition of $\Z_2\otimes \Gamma $ is 
$f_0 \perp 2 f_1 \perp 4 f_2 $ such that $f_0$ and $f_2$ are even 
unimodular forms. 
Without loss of generality we assume that $b\leq 7$ (otherwise we may
interchange the roles of $\Lambda $ and $\Gamma $).
In total this leaves 176 possibilities for the genus of $\Gamma $. 
We construct one lattice from each of the 176 genera.
Most of them were obtained as random sublattices of suitable index
of the maximal even
lattices given in Section \ref{maxeven}, representatives for 
the last 20 genera had to be constructed as sublattices of 
orthogonal sums of scaled root lattices and their duals. 
The mass of most of the genera is too big to enumerate all lattices 
in them. 
For each genus we calculate the theta series $T$ of the chosen
representative. The determinant of $\Gamma $ is 
either a square or 3 times a square since the 
dimensions of both forms ($f_0$ and $f_2$) in the
Jordan decomposition above is even and so is $\dim (f_1)$ and hence 
$a$. 
 Then we know that  $\theta _{\Gamma } = T + f $ for some
$f \in {\cal S}_{\epsilon }$, where for $\epsilon = -1$ or $\epsilon = -3$ 
the space 
${\cal S}_{\epsilon } $ is the cuspidal subspace of 
${\cal M}_{7} (\Gamma _{0}(12), \left( \frac{\epsilon }{\cdot } \right) ) $.
The condition that 
$$\begin{array}{lll} 
\theta _{\Gamma } =  & 1 + 0q^2 + 0 q^4 + 0 q^6 + 2\cdot 21\cdot 24 q^8 + \ldots  &
\mbox{ and  }  \\
\theta _{\sqrt{12} \Lambda ^* } = W_{12} (\theta _{\Gamma }) 
 =  & 1 + 0q^2 + 0 q^4 + 0 q^6 + 2\cdot 21\cdot 24 q^8 + \ldots 
\end{array}
$$
has a 1-parametric space of solutions if $\epsilon = -3$ and 
a 2-parametric space of solutions if $\epsilon = -1$. 
Using the package \cite{LRS} to search for a solution in 
this space of which the first 40 coefficients are nonnegative, 
we find no such solution in all 176 cases.
This shows that there is no suitable theta series and hence no 
such lattice $\Gamma $.
\eb

\begin{lemma}
If $s_1=6$ then $t_1 \neq 32$. 
\end{lemma}

\bew
Assume that $s_1=6$ and $t_1=32$ and rescale such that 
$\min(\Gamma ) = 4$, $\Gamma = \Lambda ^*$ with $\min (\Lambda ) = \frac{4}{3}$.
Then $\Gamma $ is even and so is $12\Lambda $. 
We proceed as in the proof of Lemma \ref{beispmodular}. 
From Lemma \ref{boundsfromgamma} we find that 
$$ 166 \leq \det(\Gamma ) = 2^a 3^b \leq 28756 .$$
Moreover $1\leq a\leq 27$ and $1\leq b \leq 13$ and the 
 Jordan decomposition of $\Z_2\otimes \Gamma $ is 
$f_0 \perp 2 f_1 \perp 4 f_2 $ such that $f_0$ and $f_2$ are even 
unimodular forms. 
In total this leaves 149 possibilities for the genus of $\Gamma $. 
We construct one lattice from each of the 149 genera.
For each genus we calculate the theta series $T$ of the chosen
representative. Again the determinant of $\Gamma $ is 
either a square or 3 times a square
and   $\theta _{\Gamma } = T + f $ for some
$f \in {\cal S}_{\epsilon }$, where for $\epsilon = -1$ or $\epsilon = -3$ 
the space 
${\cal S}_{\epsilon } $ is the cuspidal subspace of 
${\cal M}_{7} (\Gamma _{0}(12), \left( \frac{\epsilon }{\cdot } \right) ) $.
The condition that 
$$\begin{array}{ll} 
\theta _{\Gamma } =  & 1 + 0q^2 + 2\cdot 21\cdot 32 q^4 + \ldots \hfill{ 
\mbox{ and  }}  \\
\theta _{\sqrt{12} \Lambda ^* } = W_{12} (\theta _{\Gamma }) 
 =  & 1 + 0q^2 + 0 q^4 + 0 q^6 + 0 q^8 + 0 q^{10} + 0 q^{12} + 0 q^{14} + 2\cdot 21\cdot 6 q^{16} + \ldots 
\end{array}
$$
determines $\theta _{\Gamma }$ uniquely. 
One finds no solution in all cases where 
 $\epsilon = -3$ and a unique solution 
 if $\epsilon = -1$. 
This unique modular form has negative coefficients in all cases.
This shows that there is no suitable theta series and hence no 
such lattice $\Gamma $.
\eb

\renewcommand{\arraystretch}{1}
\renewcommand{\baselinestretch}{1}
\large


\normalsize

\begin{thebibliography}{10}
\bibitem{Magma} J. Cannon,  MAGMA Computational Algebra System  \\
http://magma.maths.usyd.edu.au/magma/htmlhelp/MAGMA.htm
\bibitem{Sage} W. Stein, Sage, Open source mathematics software, \\
http://modular.math.washington.edu/sage/
\bibitem{LRS} D. Avis, lrs Vertex Enumeration/Convex Hull package \\
http://cgm.cs.mcgill.ca/$\sim$avis/C/lrs.html
\bibitem{Anstreicher} 
K. M. Anstreicher, {\it Improved linear programming bounds for antipodal spherical codes.} Discrete Comput. Geom. 28 (2002), no. 1, 107--114
\bibitem{BachocVenkov} C. Bachoc, 
B. Venkov, {\it Modular forms, lattices and spherical designs.}
Monogr. Ens. Math. vol. 37, (2001) 87-111 
\bibitem{tight} 
Eiichi Bannai, A. Munemasa, B. Venkov,
 {\em The nonexistence of certain tight spherical designs.}
Algebra i Analiz  16  (2004),  no. 4, 1--23;  translation in 
 St. Petersburg Math. J.  16  (2005),  no. 4, 609--625 
\bibitem{Elkies}  H. Cohn, N.D. Elkies,
{\it New upper bounds on sphere packings I}, Annals of Math. 157 (2003), 689-714 (math.MG/0110009). 
\bibitem{mass} J.H. Conway, N.J A. Sloane, {\it 
Low-Dimensional Lattices IV: The Mass Formula},  Proc. Royal Soc. London, Series A, 419 (1988), pp. 259-286
\bibitem{SPLAG} J.H. Conway, N.J A. Sloane, {\it Sphere Packings, Lattices
and Groups.} 3rd edition, Springer-Verlag 1998.
\bibitem{coulange1}  R. Coulangeon,
{\it Spherical designs and zeta functions of lattices}, International Mathematics Research Notices (2006).
\bibitem{Coulange2}  R. Coulangeon,
{\it On Epstein zeta function of Humbert forms},
 International Journal of Number Theory 4 (2008), no. 3, 387-401.
\bibitem{Delone} B.N. Delone, S.S. Ryskov,
{\it A contribution to the theory of the extrema of a multidimensional 
$\zeta $-function.} Doklady Akademii Nauk SSSR {\bf 173} (1967), 991-994.
translated as Soviet Math. Doklady {\bf 8}  (1967), 499-503.  
\bibitem{king} O. King, {\it A mass-formula for unimodular
lattices with no roots.} Math. Comp. 72 (2003) 839-863
\bibitem{Kneser}
M. Kneser,
{\it Klassenzahlen definiter quadratischer Formen},
Archiv der Math. 8 (1957) 241--250.
%\bibitem{Levenstein} Levenstein bound 
\bibitem{Martinet} J. Martinet, {\it
Les R\'eseaux parfaits des espaces Euclidiens.}
Masson (1996)
\bibitem{Nebehome} % Gabriele Nebe's homepage 
http://www.math.rwth-aachen.de/homes/Gabriele.Nebe/
\bibitem{dim10} G. Nebe, B. Venkov, {\it The strongly perfect lattices of dimension 10.}
J. Th\'eorie de Nombres de Bordeaux {12} (2000) 503-518.
\bibitem{dim12} G. Nebe, B. Venkov, {\it Low dimensional strongly perfect lattices I: The 12-dimensional case.}
 L'enseignement Math\'ematique, 51 (2005) 129-163
\bibitem{dim13} G. Nebe, B. Venkov, {\it Low dimensional strongly perfect lattices II: The 13-dimensional case.} (in preparation)
\bibitem{ScharlauHemkemeier} R. Scharlau, B. Hemkemeier, 
{\it Classification of integral lattices with large class number.} Math. Comput. 67(222): 737-749 (1998)
\bibitem{Scharlau} W. Scharlau, {\it Quadratic and Hermitian forms.}
Springer Grundlehren {\bf 270} (1985)
\bibitem{Schuermann} A. Sch\"urmann, {\it Perfect, strongly eutactic 
lattices are periodic extreme.} (preprint 2008)
\bibitem{Sobolev} S.L. Sobolev, {\it Formulas for mechanical cubatures in $n$-dimensional space.} Doklady Akademii Nauk SSSR {\bf 137} (1961), 527-530. 
\bibitem{Venkov} B. Venkov, {\it R\'eseaux et designs sph\'eriques.}
Monogr. Ens. Math. vol. 37, (2001) 10-86
\end{thebibliography}
\end{document}